\documentclass{article}

\usepackage{mathmag}

\usepackage{amsmath,amsthm}     
\usepackage{graphicx}     
\usepackage{hyperref} 
\usepackage{url}
\usepackage{amsfonts} 

\usepackage{multicol}

\theoremstyle{theorem}

\theoremstyle{definition}

\allowdisplaybreaks

\makeatletter
\@addtoreset{footnote}{page}
\makeatother

\usepackage{enumerate}
\usepackage{subcaption}

\begin{document}

\newpage
\title{Fault-Free Tileability of Rectangles, Cylinders, Tori, and M\"{o}bius Strips with Dominoes}

\author{Emily Montelius\\               
\scriptsize Department of Mathematical Sciences, Coe College\\    
1220 1$^{\text{st}}$ Ave NE, Cedar Rapids, Iowa\\                
                 }    

\maketitle

\section{Introduction}

\begin{figure}[h]
     \centering
     \begin{subfigure}[b]{0.3\textwidth}
         \centering
         \includegraphics[width=\textwidth]{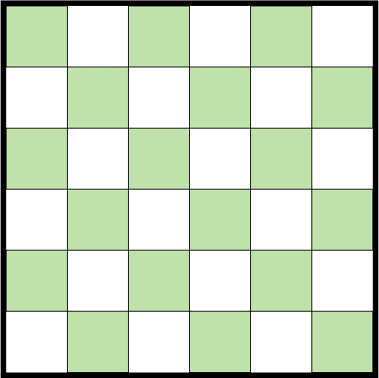}
         \caption{}
         \label{fig:6x6blank}
     \end{subfigure}
     \hspace*{1cm}
     \begin{subfigure}[b]{0.3\textwidth}
         \centering
         \includegraphics[width=\textwidth]{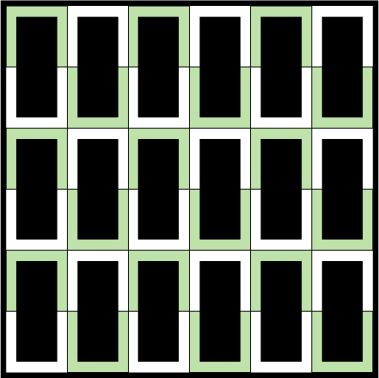}
         \caption{}
         \label{fig:6x6tiled}
     \end{subfigure}
     
        \caption{A $6\times 6$ board and one tiling of the board}
        \label{fig:6x6examples}
\end{figure}

We begin with a \textit{board}, or an $a\times b$ rectangular grid, and we place \textit{tiles} on the board. We use $1\times 2$ dominoes as our tiles. Each $1\times 1$ unit of area on the board is a \textit{space}. We say a board is \textit{tileable} if it can be filled with tiles such that no space is left uncovered, nor any space covered by more than $1$ tile. For example, the $6\times 6$ board is tileable, as shown in Figure~\ref{fig:6x6examples}. Not all boards are tileable by dominoes. When a tiling is possible, a working solution must be shown, otherwise a proof must be given for why the board is not tileable. For example, the $5\times5$ board is not tileable. Try to draw tiles in Figure~\ref{fig:5x5}, and formulate an explanation as to why the $5\times5$ is not tileable. An explanation is given below.

The $5\times 5$ board has $25$ spaces to be covered. Each domino covers two spaces. Thus, with dominoes, we can only cover an even number of spaces and cannot tile a board with an odd area. This method of proof, known as \textit{parity}, causes all odd by odd boards to be not tileable. 

\begin{figure}[h]
     \centering
         \includegraphics[width=0.8\textwidth]{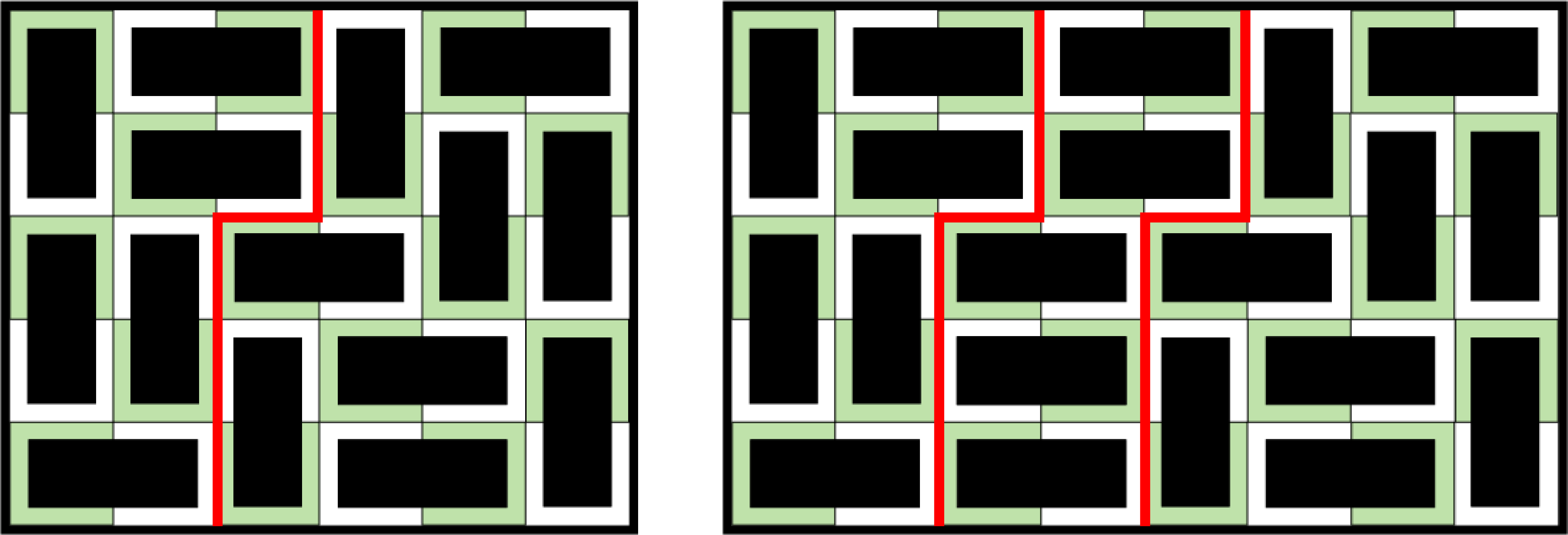}
         \caption{A $5\times 6$ board and one fault-free tiling of the board with an expansion}
         \label{fig:expansion}
\end{figure}

Assuming tiles extend to the edges of the spaces, a \textit{fault-line} exists anywhere a board could be folded without the crease crossing through the interior of a tile. We call a tiling with no fault-lines \textit{fault-free}. We call a board for which a fault-free tiling exists \textit{fault-free tileable}. In Figure~\ref{fig:6x6examples} we see a board with $7$ fault-lines, and in Figure~\ref{fig:expansion} we see $2$ fault-free tileable boards. Our goal is to determine which boards are or are not fault-free tileable.

Previous work on fault-free tiling of rectangular boards with dominoes was completed by Ron Graham in ``Fault-Free Tilings of Rectangles"~\cite{Graham}. Here we identify a few techniques and examples used by Graham to determine whether a board is fault-free tileable or not. The first example we will look at is how known fault-free boards can be expanded. 

We begin with the $5\times 6$ board, which is fault-free tileable as shown in Figure~\ref{fig:expansion}. It can be split into two non-rectangular boards across the bold line and then expanded by multiples of $2$. By placing the new tiles in this way, parallel, lined up, and offset by one, all fault-lines will remain covered. We can expand the board, in this way, arbitrarily many times, and in either direction. In general, if we know an $a\times b$ board is fault-free tileable, then we know the board of size $(a+2n)\times (b+2m)$, where $n,m \in \mathbb{N}$, is also fault-free tileable.

\begin{figure}[h]
\centering
         \subcaptionbox{\label{fig:6x6ex1}}{\includegraphics[scale=1.15]{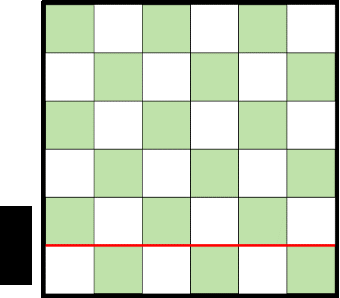}}\hspace*{1cm}
         \subcaptionbox{\label{fig:6x6ex2}}{\includegraphics[scale=.25]{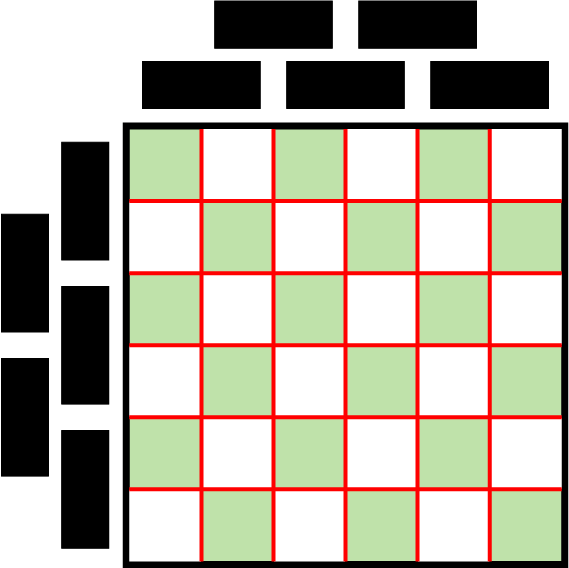}}
         
         \vspace*{.5cm}
         \subcaptionbox{\label{fig:6x6ex3}}{\includegraphics[scale=.25]{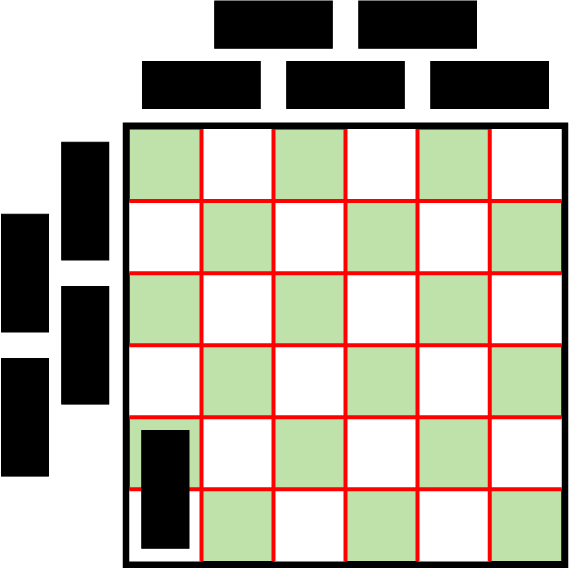}}\hspace*{1cm}
         \subcaptionbox{\label{fig:6x6ex4}}{\includegraphics[scale=.25]{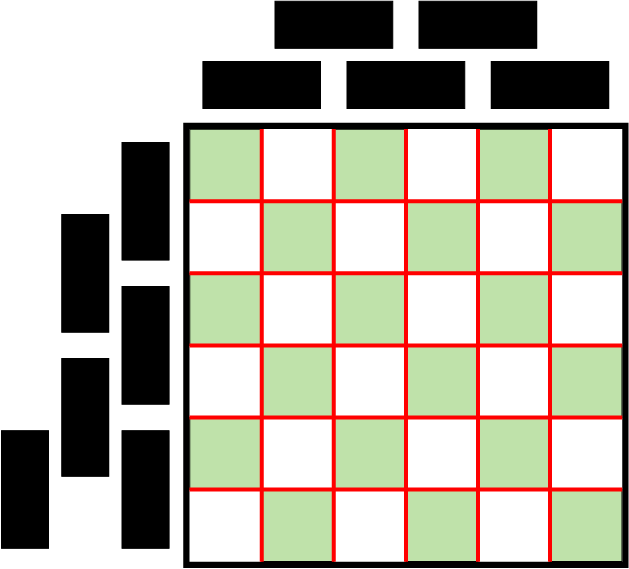}}
         
         \vspace*{.5cm}
         \subcaptionbox{\label{fig:6x6ex5}}{\includegraphics[scale=.25]{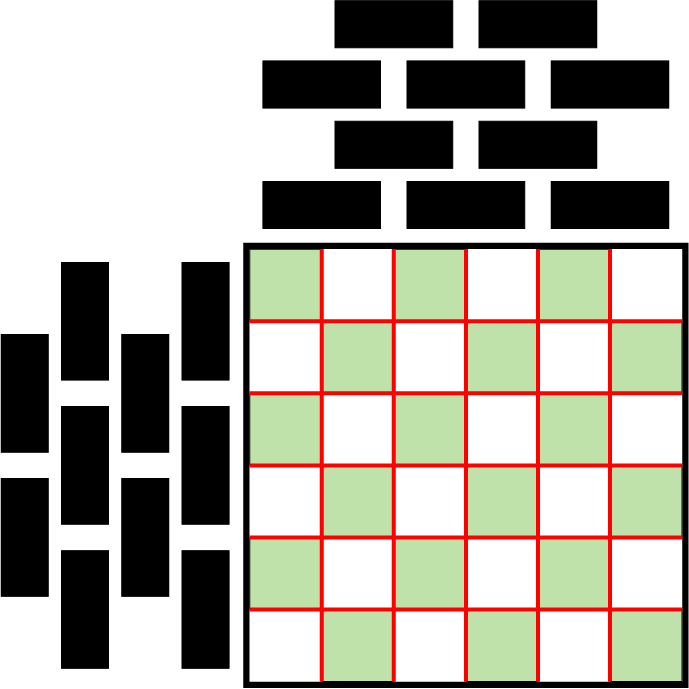}}

        \caption{The $6\times 6$ board and its required tiles to be fault-free tileable}
        \label{fig:6x6fftexamples}
\end{figure}

Another example Graham provides is a proof that the $6\times 6$ board, a board we know to be tileable, is not fault-free tileable. To determine why, we start by noting fault-lines that must be covered. At least $1$ tile must cross each fault-line, so we can place tiles to the side of the board, which are fixed to move along the potential fault-line. This is shown for one fault-line in Figure~\ref{fig:6x6ex1}. It is required this tile be placed somewhere on it's corresponding fault-line in order to create a fault-free tiling of the board. In Figure~\ref{fig:6x6ex2}, we have $10$ tiles placed because there are 10 fault lines on the board. Looking at any edge-most column or row, we can see there are $6$ spaces in that column or row that must be covered by a tile. For example, let's consider the bottom row. The $6$ spaces in the row will be covered in either pairs or individually, depending on the orientation of the domino. A horizontal tile would cover $2$ spaces in the row, while a vertical tile would cover one space. As shown in Figure~\ref{fig:6x6ex3}, with the currently required tiles, there is one individually covered space in the bottom row. It is covered by the vertical domino required to cover the bottom most fault-line. This leaves $5$ spaces left to be covered. Because $5$ is an odd number, these spaces cannot be covered only by horizontal tiles. Thus, one or more tiles of vertical orientation must be added to the bottom row and will also cross the bottom fault-line. For this reason, we add one additional tile to the side of the board, as in Figure~\ref{fig:6x6ex4}. This pattern continues for each row and column with an odd number of spaces remaining to be covered, leading to $20$ tiles required to cover all fault-lines, as shown in Figure~\ref{fig:6x6ex5}. However, the $6\times 6$ board only has room for $18$ tiles, because there are $36$ spaces, showing it is impossible to fault-free tile this board because more tiles are required than spaces allow for. This investigation of the number of required tiles for a fault-free tiling is an important concept that will be used frequently in other cases to  determine boards as not fault-free tileable.

\section{Cylinders}

\begin{figure}[h]
     \centering
         \includegraphics[scale=.25]{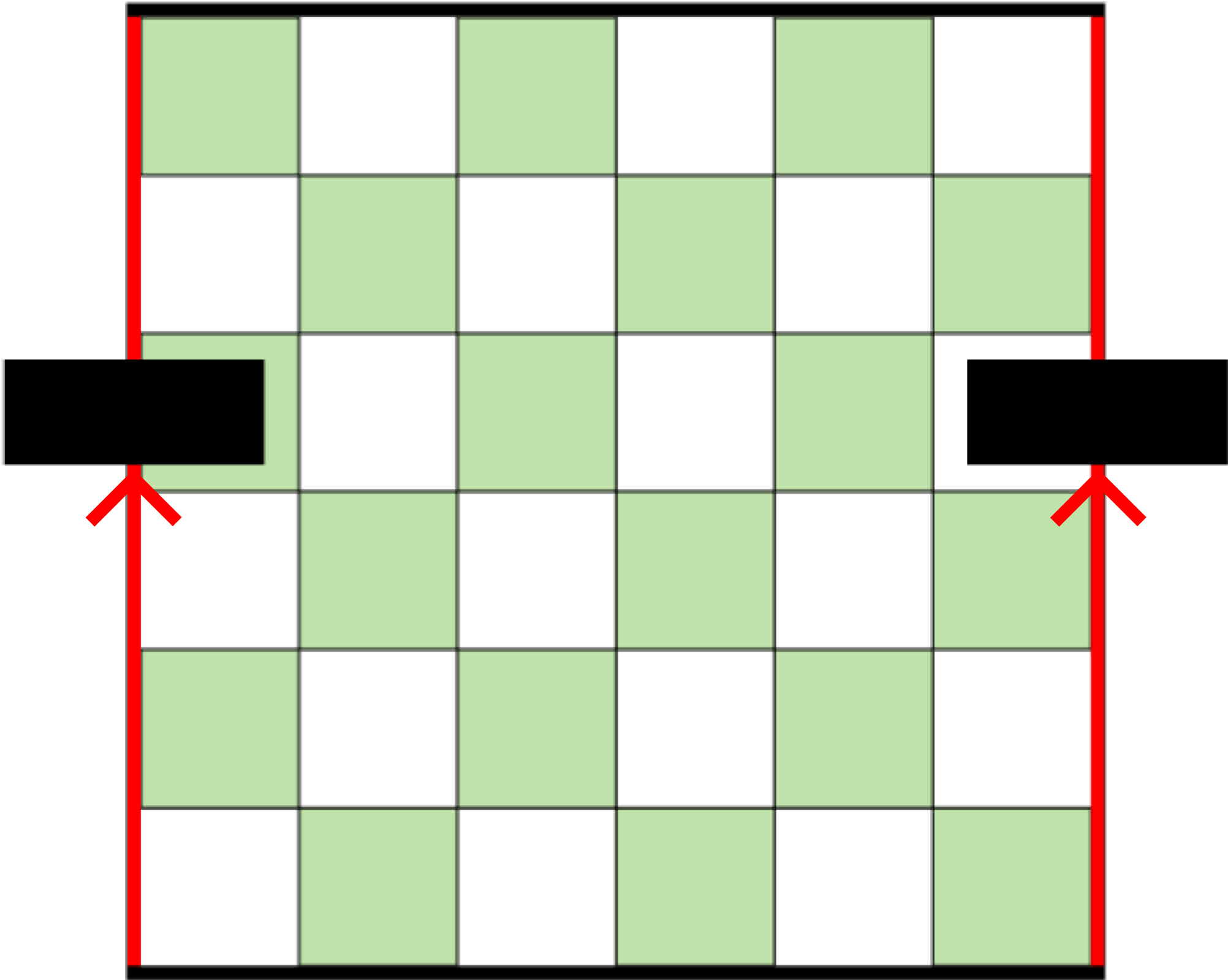}
         \caption{A $6^\prime \times 6$ cylindrical board with one tile placed}
         \label{fig:cylinderex}
\end{figure}

\begin{figure}[h]
     \centering
         \includegraphics[width=0.3\textwidth]{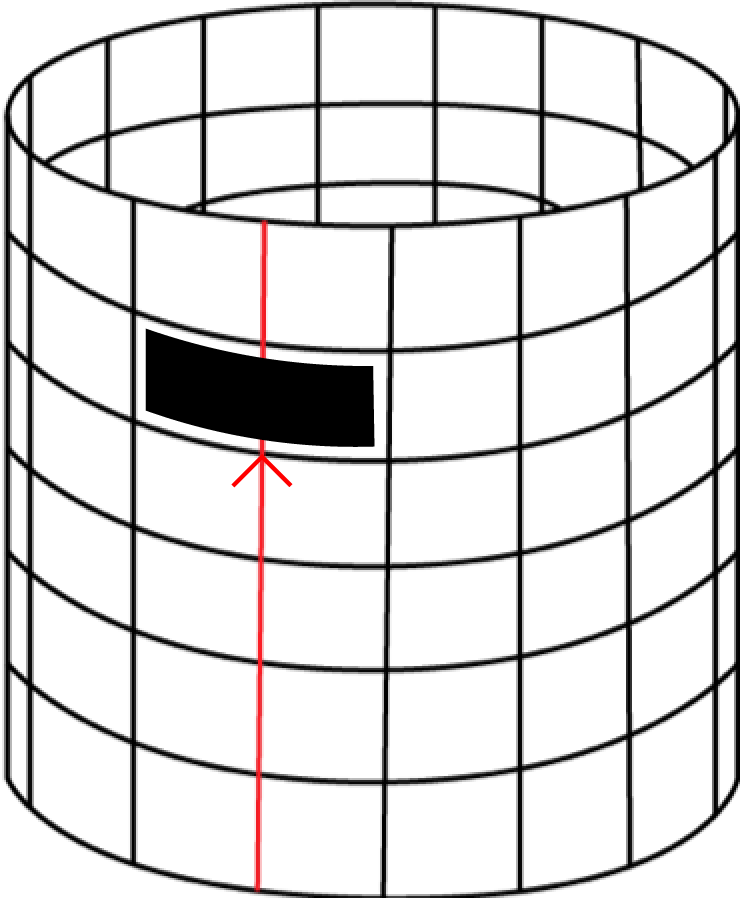}
         \caption{A $6^\prime \times 15$ cylinder board shown in three dimensions}
         \label{fig:literalcylinder}
\end{figure}

We now proceed to the original part of this work, determining which shapes other than rectangles are fault-free tileable. We begin by investigating cylinders. We imagine a cylindrical board as a rectangular grid where one set of opposing edges are glued together. We show these cylinders in two dimensions where tiles wrap around the sides denoted with arrows, as shown in Figure~\ref{fig:cylinderex}. Though it appears there are two tiles in Figure~\ref{fig:cylinderex}, this is in fact only one tile. The vertical edges are actually one fault-line that appears on both sides of the flattened cylinder's image. A cylinder with one tile and a vertical fault line denoted with an arrow is shown in three dimensions in Figure~\ref{fig:literalcylinder}. We will denote a cylindrical board as $a^\prime \times b$ where $a$ is the height and $b$ is the circumference of the cylinder. For generality, $a$ is the length of vertical edges, while $b$ is the length of the horizontal edges when shown in two dimensions. The prime indicates that the sides of that length are not edges, but in fact a single fault-line appearing on both sides, and this is where the image wraps around. Compared to an $a\times b$ rectangular board, the cylindrical board, $a^\prime \times b$, will have one more fault-line than the $a\times b$ board because the two vertical edges become an additional fault-line.

\begin{figure}[h]
\centering
         \subcaptionbox{\label{fig:6x4'}}{\includegraphics[scale=.25]{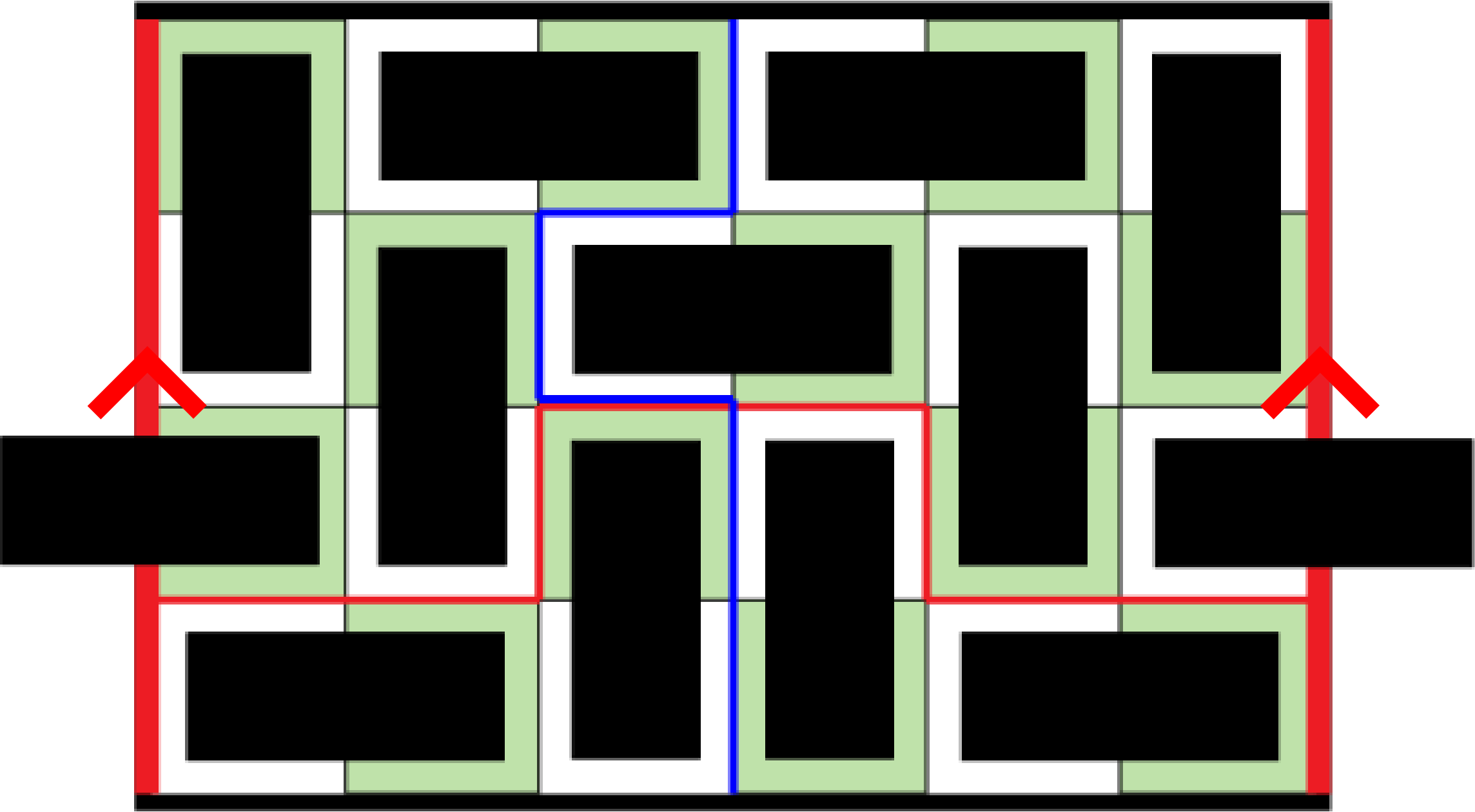}}
         
         \vspace*{.5cm}
         \subcaptionbox{\label{fig:6x7'}}{\includegraphics[scale=.25]{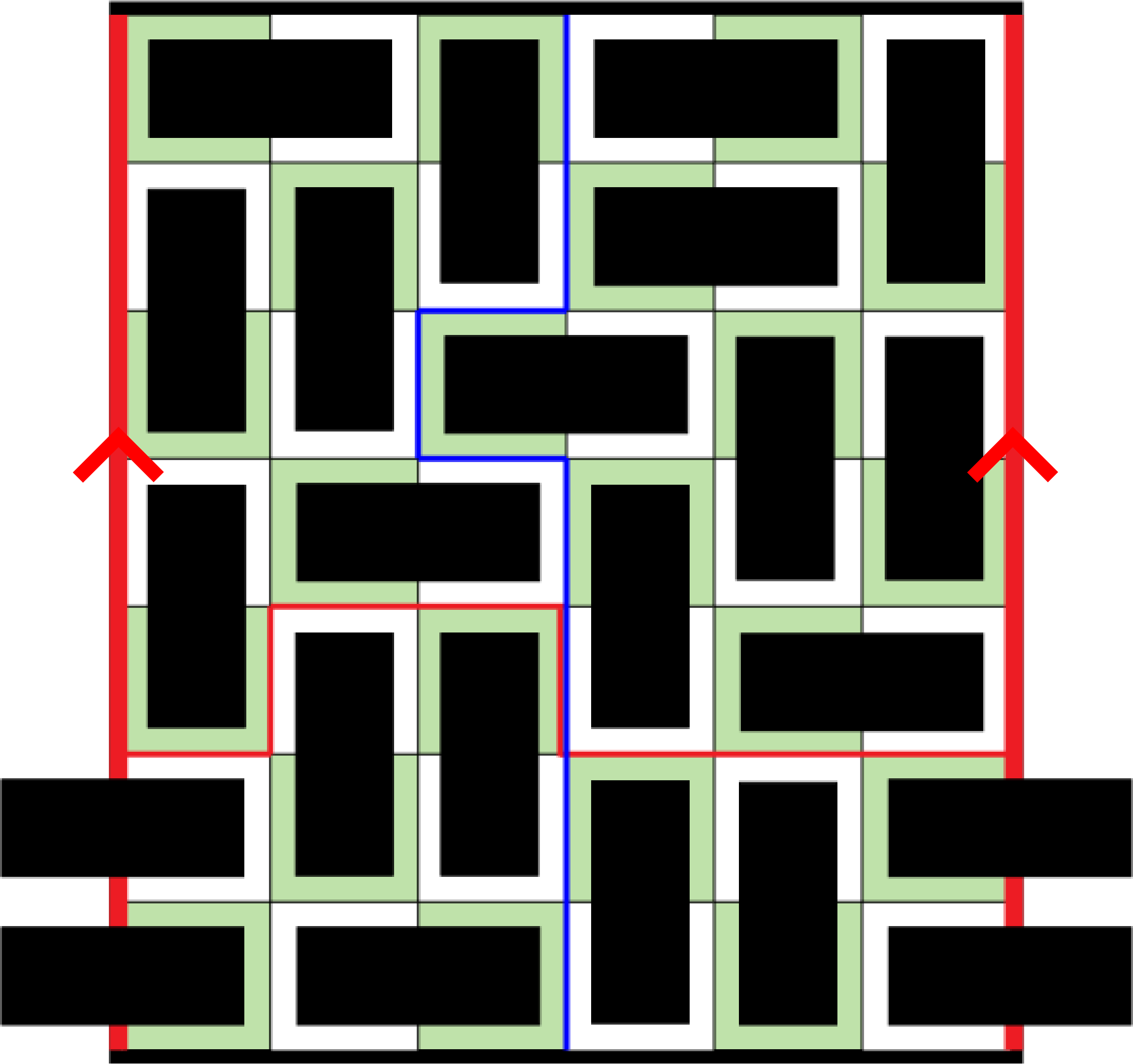}}\hspace*{1cm}
          \subcaptionbox{\label{fig:7x6'}}{\includegraphics[scale=.25]{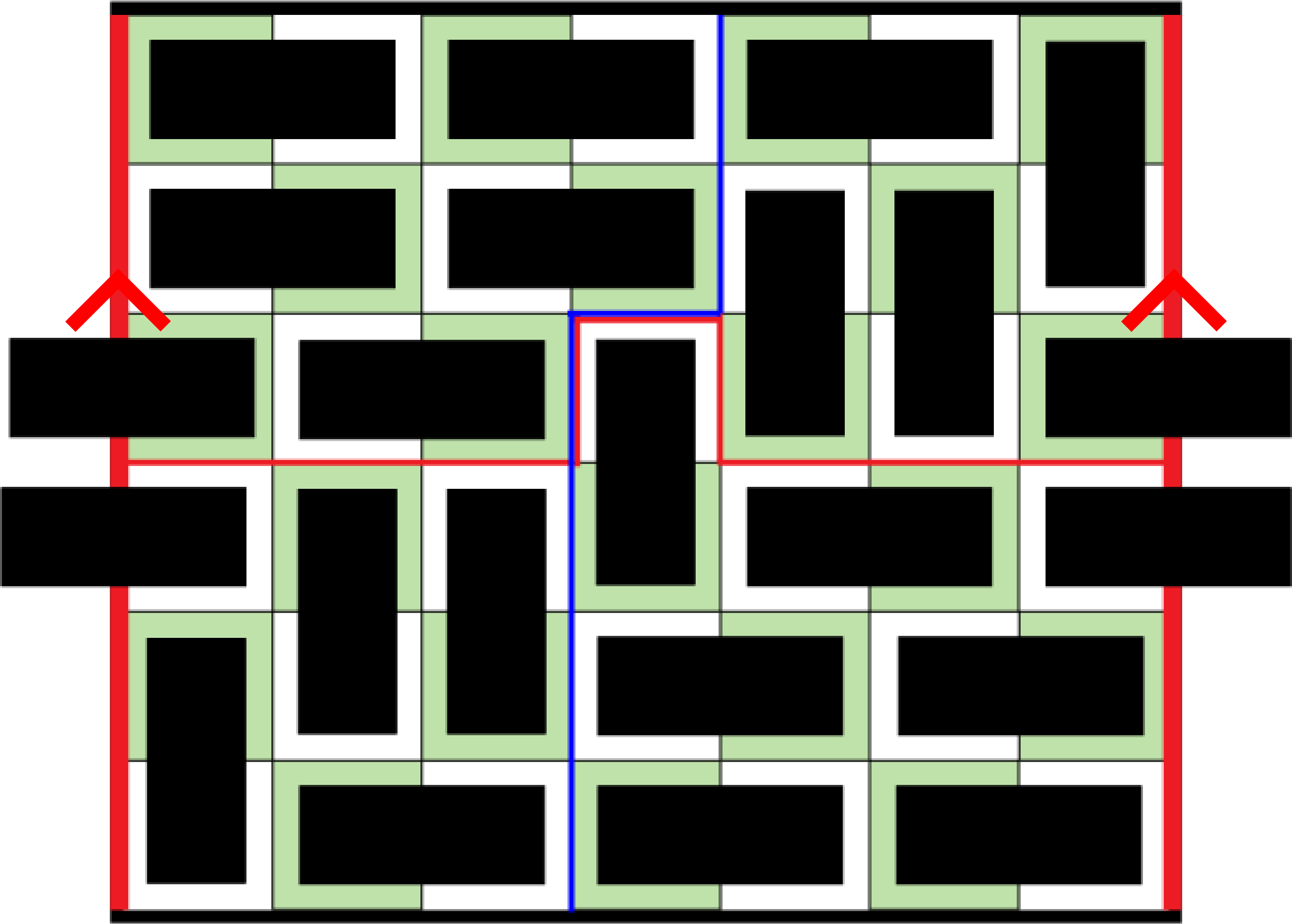}}
          
         \vspace*{.5cm}
          \subcaptionbox{\label{fig:5x8'}}{\includegraphics[scale=.25]{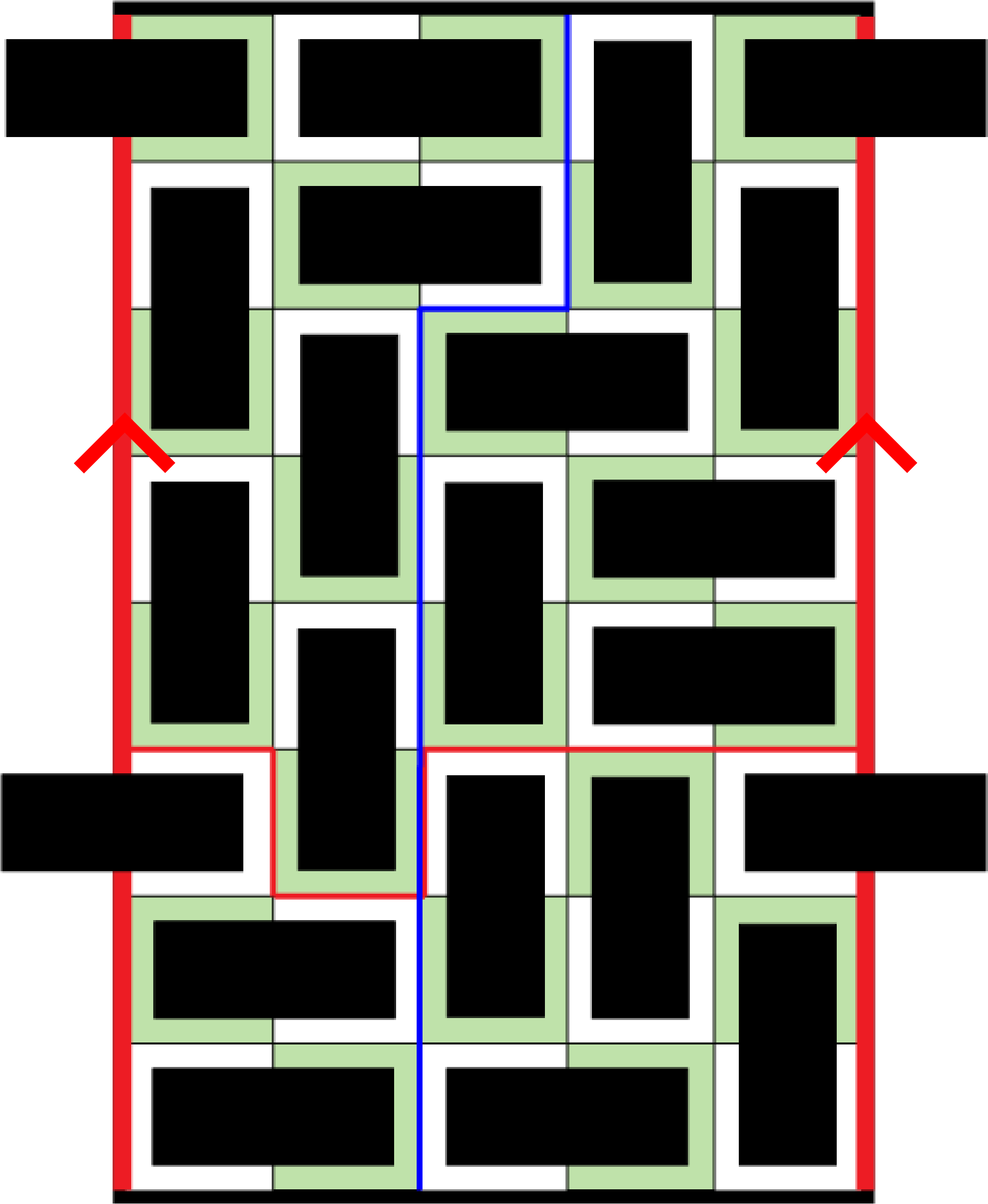}}\hspace*{1cm}
         \subcaptionbox{\label{fig:8x5'}}{\includegraphics[scale=.25]{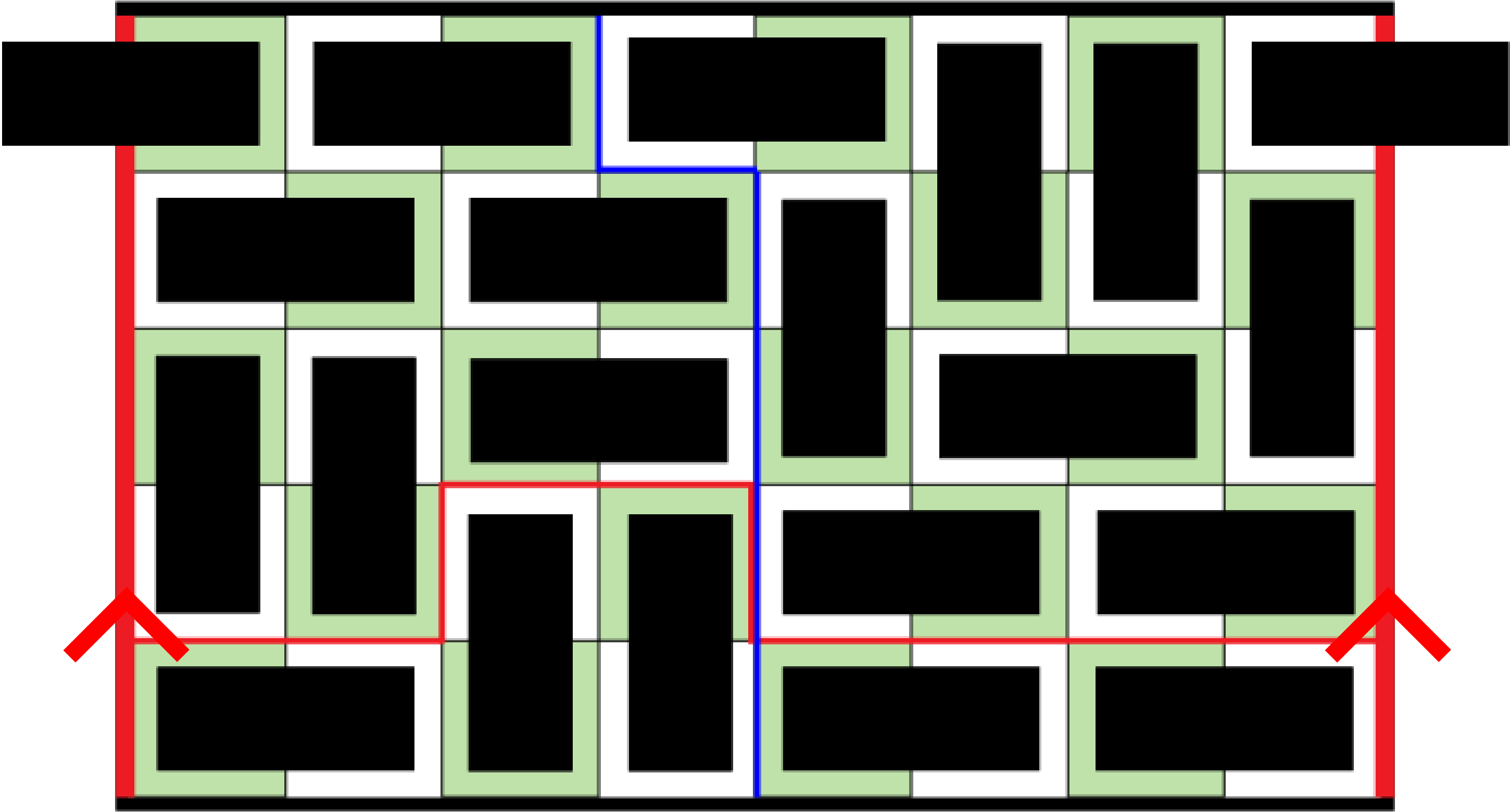}}

        \caption{The fault-free tileable cylindrical boards}
        \label{fig:tileablecyl}
\end{figure}

The following cylinders are fault-free tileable as shown in Figure~\ref{fig:tileablecyl}:
\begin{enumerate}[(a)]
    \item $4^\prime \times 6$
    \item $7^\prime \times 6$
    \item $6^\prime \times 7$
    \item $8^\prime \times 5$
    \item $5^\prime \times 8$
\end{enumerate}

Again, these can be expanded by multiples of two in either direction, just as the rectangular boards could be. Thus the following cylinders are fault-free tileable, where $n,m \in \mathbb{N}$:
\begin{enumerate}[(a)]
    \item $(4+2n)^\prime \times (6+2m)$
    \item $(7+2n)^\prime \times (6+2m)$
    \item $(6+2n)^\prime \times (7+2m)$
    \item $(8+2n)^\prime \times (5+2m)$
    \item $(5+2n)^\prime \times (8+2m)$
\end{enumerate}

We now show why the following boards are not fault-free tileable.
\begin{itemize}
    \item $\text{odd}^\prime \times \text{odd}$
    \item $(2n)^\prime \times 1$
    \item $1^\prime \times (2n)$
    \item $(2+n)^\prime \times 2$
    \item $2^\prime \times (2+n)$
    \item $(4+2n)^\prime \times 3$
    \item $3^\prime \times (4+2n)$
    \item $(4+n)^\prime \times 4$
    \item $4^\prime \times (5+2n)$
    \item $6^\prime \times 5$
    \item $5^\prime \times 6$
\end{itemize}

The odd by odd boards are not fault-free tileable by parity, just as the rectangular boards could not be.

The $(2n)^\prime \times 1$ and $1^\prime \times (2n)$ boards are not fault-free tileable, as any attempt to place a first tile will result in a fault-line.

The $(2+n)^\prime \times 2$ boards are not fault-free tileable. All vertical fault-lines must be crossed. Arbitrarily picking one of the two vertical fault lines, a horizontal tile must be placed to cover it. This will result in a horizontal fault-line either above, below, or above and below the tile once it is placed.

Similarly, the $2^\prime \times (2+n)$ boards are not fault-free tileable. The horizontal fault-line must be crossed. Placing a vertical tile to cover this will result in a vertical fault-line on either side of the tile.

\begin{figure}[h]
     \centering
         \includegraphics[scale=.25]{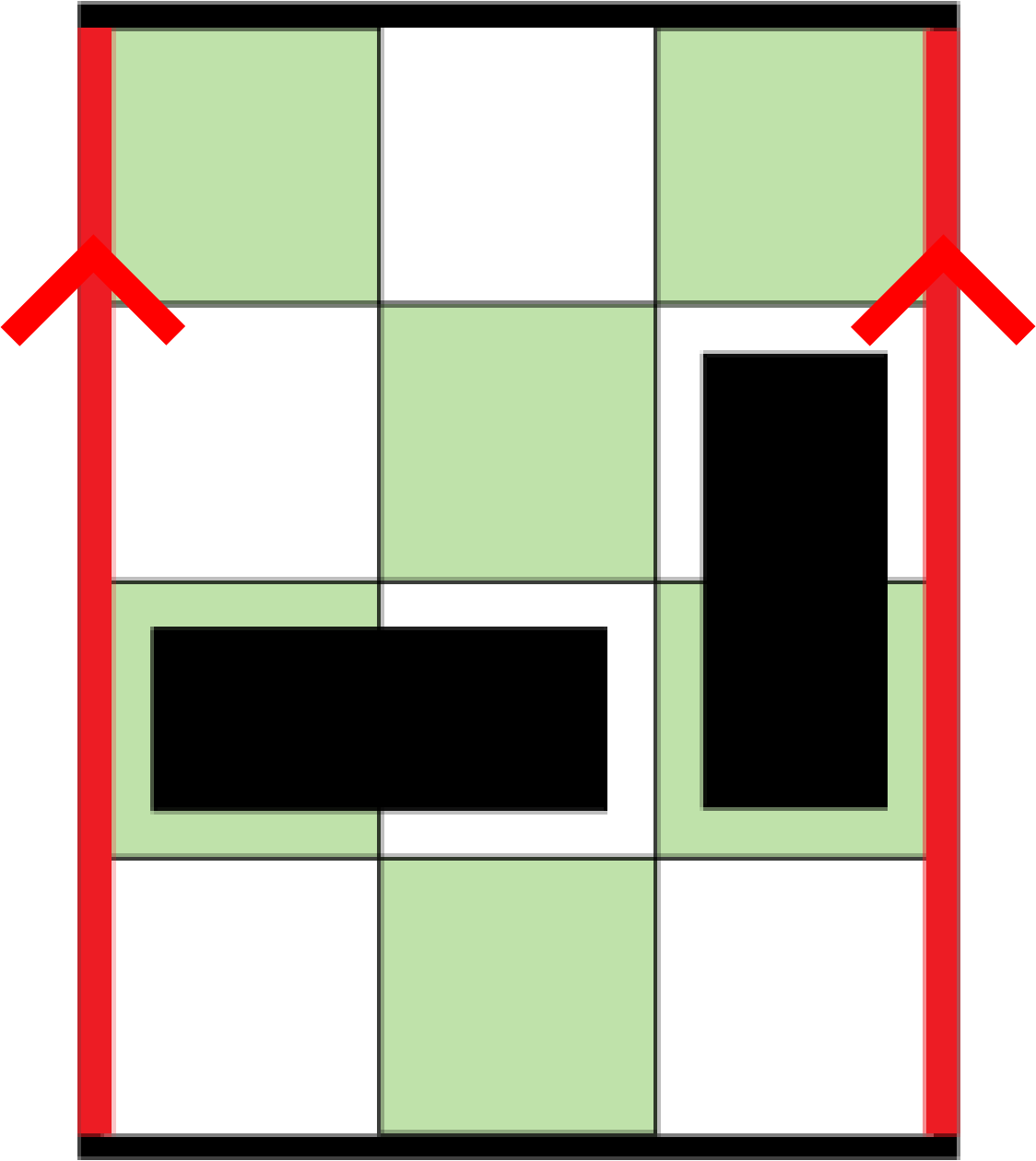}
         \caption{The $4^\prime \times 3$ board and its forced tile sequence}
         \label{fig:3x4'}
\end{figure}

The $(4+2n)^\prime \times 3$ boards are not fault-free tileable. For example, consider the $4^\prime \times 3$ cylinder shown in Figure~\ref{fig:3x4'}. Three horizontal tiles are needed to cross the the three vertical fault-lines. An edge most row cannot contain more than one horizontal tile. Thus, at least one of the three vertical fault-lines must be crossed by a horizontal tile in a central row. This placement of a horizontal tile in a central row forces a vertical tile to be placed in the remaining space of that row, which will result in a horizontal fault-line. This is the case no matter how tall the cylinder becomes beyond a height of $4$.

\begin{figure}[h]
     \centering
         \includegraphics[scale=.25]{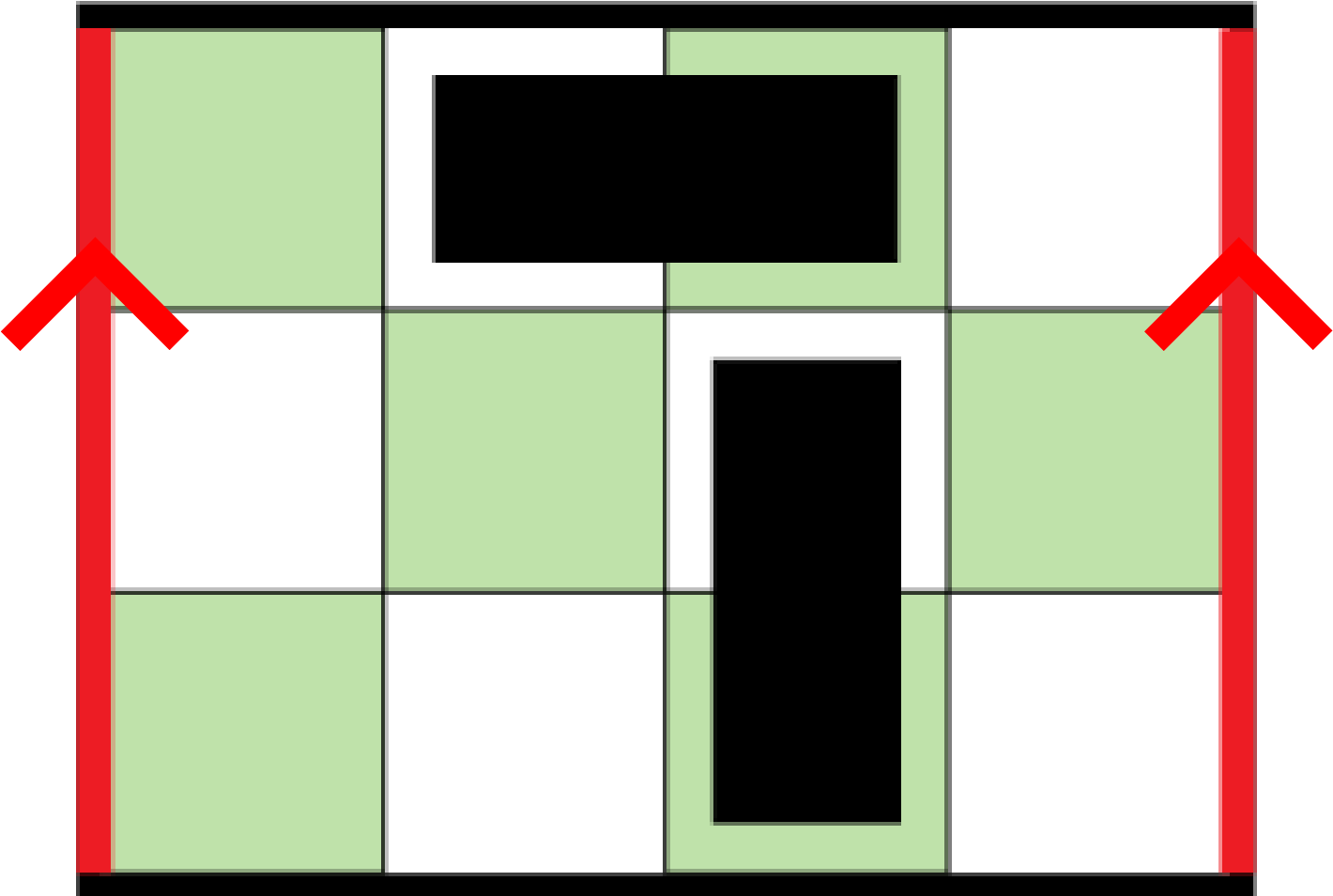}
         \caption{The $3^\prime \times 4$ board and its forced tile sequence}
         \label{fig:4x3'}
\end{figure}

The $3^\prime \times (4+2n)$ boards are not fault-free tileable. This is demonstrated in the $3^\prime \times 4$ depicted in Figure~\ref{fig:4x3'}. A horizontal fault-line must be crossed by placing a vertical tile. This will force a horizontal tile to be placed in the column to cover the one remaining space, which will result in a vertical fault-line. One may think we could place the vertical tile in an edge most column to avoid this, however this board is a cylinder and there are no vertical edges for us to place the initial tile against.

\begin{figure}[h]
     \centering
         \includegraphics[scale=.25]{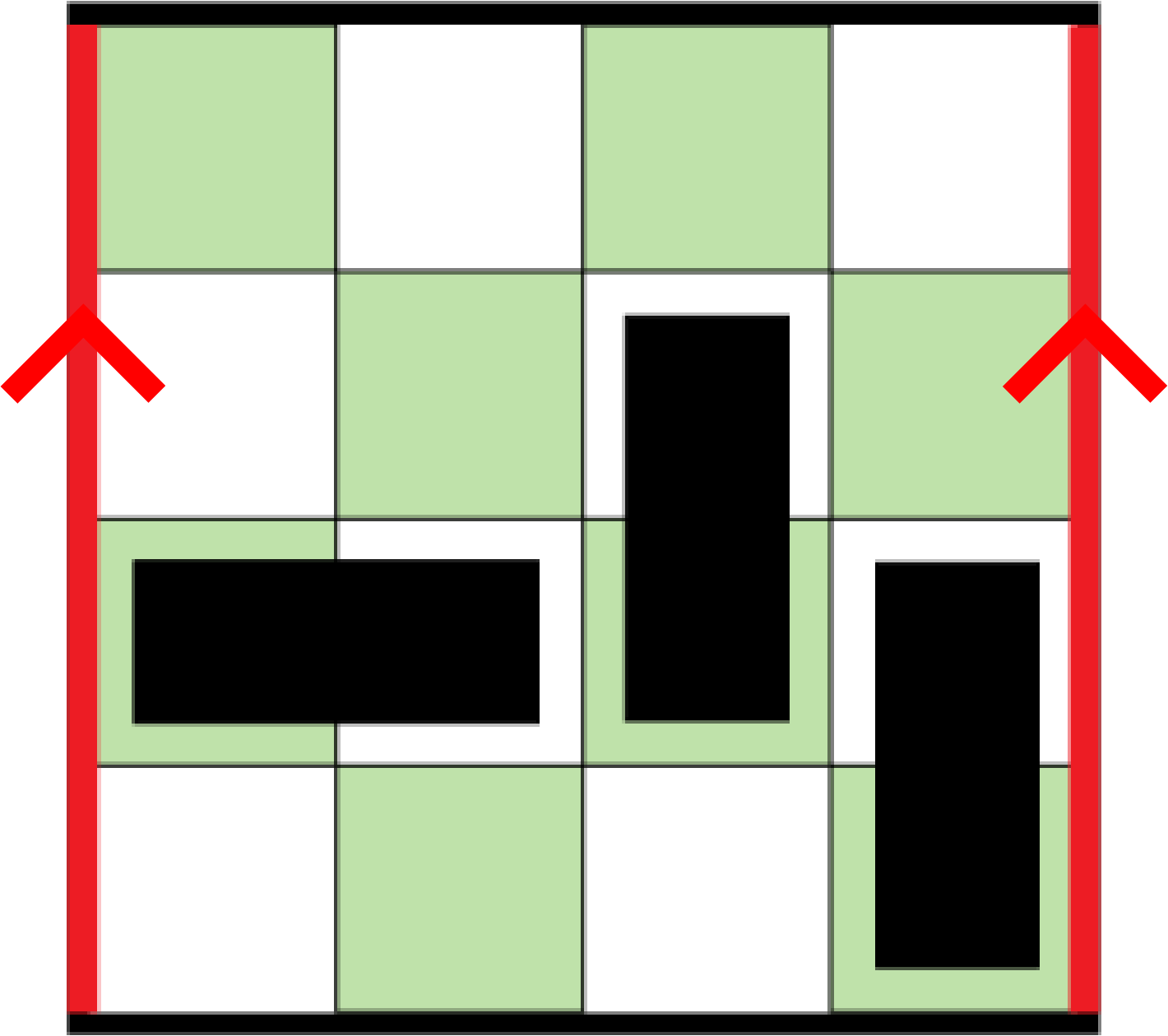}
         \caption{The $4^\prime \times 4$ board and its forced tile sequence}
         \label{fig:4x4'}
\end{figure}

The $(4+n)^\prime \times 4$ boards are not fault-free tileable. For example, Figure~\ref{fig:4x4'} shows the $4^\prime \times 4$ case. A vertical fault-line must be crossed. As many as two vertical fault-lines may be crossed at the top-most and bottom-most rows. Crossing any more than two at the edges will create a horizontal fault-line, so at least one vertical fault-line must be crossed in a central row with a horizontal tile. This will force the next two tiles to be vertical and in opposing directions to avoid creating a horizontal fault-line. This results in an odd number of spaces remaining both above and below the placed tiles. Thus, the collection of spaces created on either side are not tileable by parity. Thus the whole board is not fault-free tileable.

\begin{figure}[h]
     \centering
         \includegraphics[scale=.25]{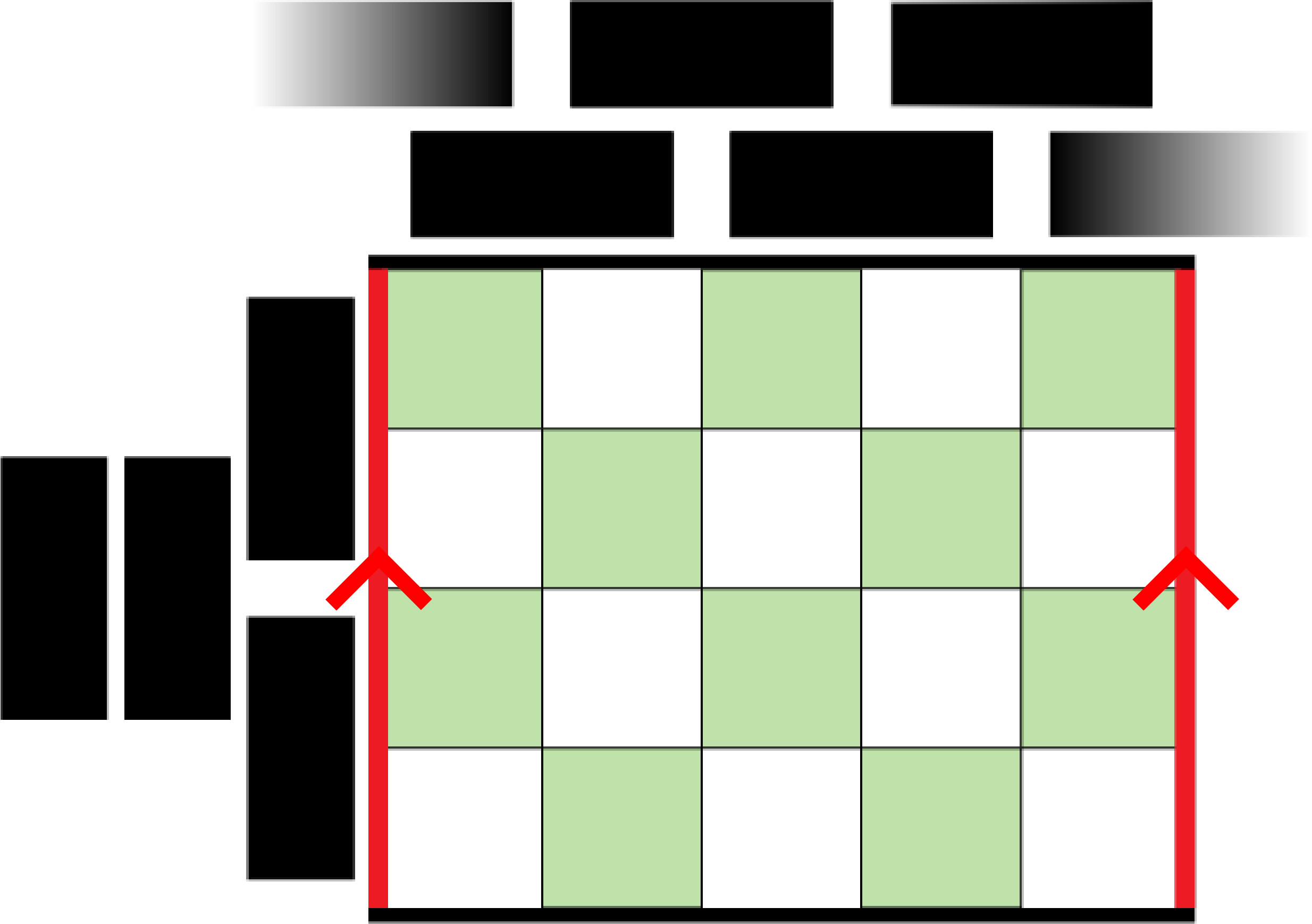}
         \caption{The $4^\prime \times 5$ board and its required tiles to be fault-free tileable}
         \label{fig:5x4'}
\end{figure}

The $4^\prime \times (5+2n)$ boards are not fault-free tileable. This is shown in a $4^\prime \times 5$ example in Figure~\ref{fig:5x4'}. To avoid accidental double-counting of tiles, any tile that appears twice has a gradient  on each copy to emphasize that it is in fact one tile. These boards require $5+2n$ tiles to cover all vertical fault-lines.  This board also requires $4$ tiles to cover the horizontal fault-lines. It requires $4$ rather than $3$ for the same reasons the $6\times 6$ required $20$ rather than $10$: the odd number of spaces remaining uncovered in the central rows. This totals to $9+2n$ tiles required, an odd number. A $ 4^\prime \times (5+2n)$ board, tiled fully, will contain $10+4n$ tiles, an even number. If we add an additional tile to the $9+2n$ required, we change the number of spaces remaining in the corresponding rows or columns by one. This means a second additional tile of the same orientation will have to be added to the same fault-line in order to maintain the affected rows' or columns' parity. For this reason we can add tiles to the board in pairs, while we cannot add a singular tile to the board. We can also add additional tiles to the board by adding a horizontal set of $5+2n$ in the same orientation as the current required tiles. Adding an additional $5+2n$ in this way is possible because it would maintain the parity of the columns. Tiles can only be added to the board in sets of $2x+(5+2n)y$, $x,y \in \mathbb{N}$. There is no other way to add tiles without creating a row or column with an odd number of remaining spaces to be covered. To get an even number of tiles to fill the board, we must add an odd number of tiles to the $9+2n$ already required. Given our restrictions on adding tiles, the only way to add an odd number is to add a new set of $5+2n$. This will give us $14+4n$ tiles in the board. This is more than the $10+4n$ there is space for. Thus it cannot be fault-free tiled.

\begin{figure}[h]
     \centering
         \includegraphics[scale=.25]{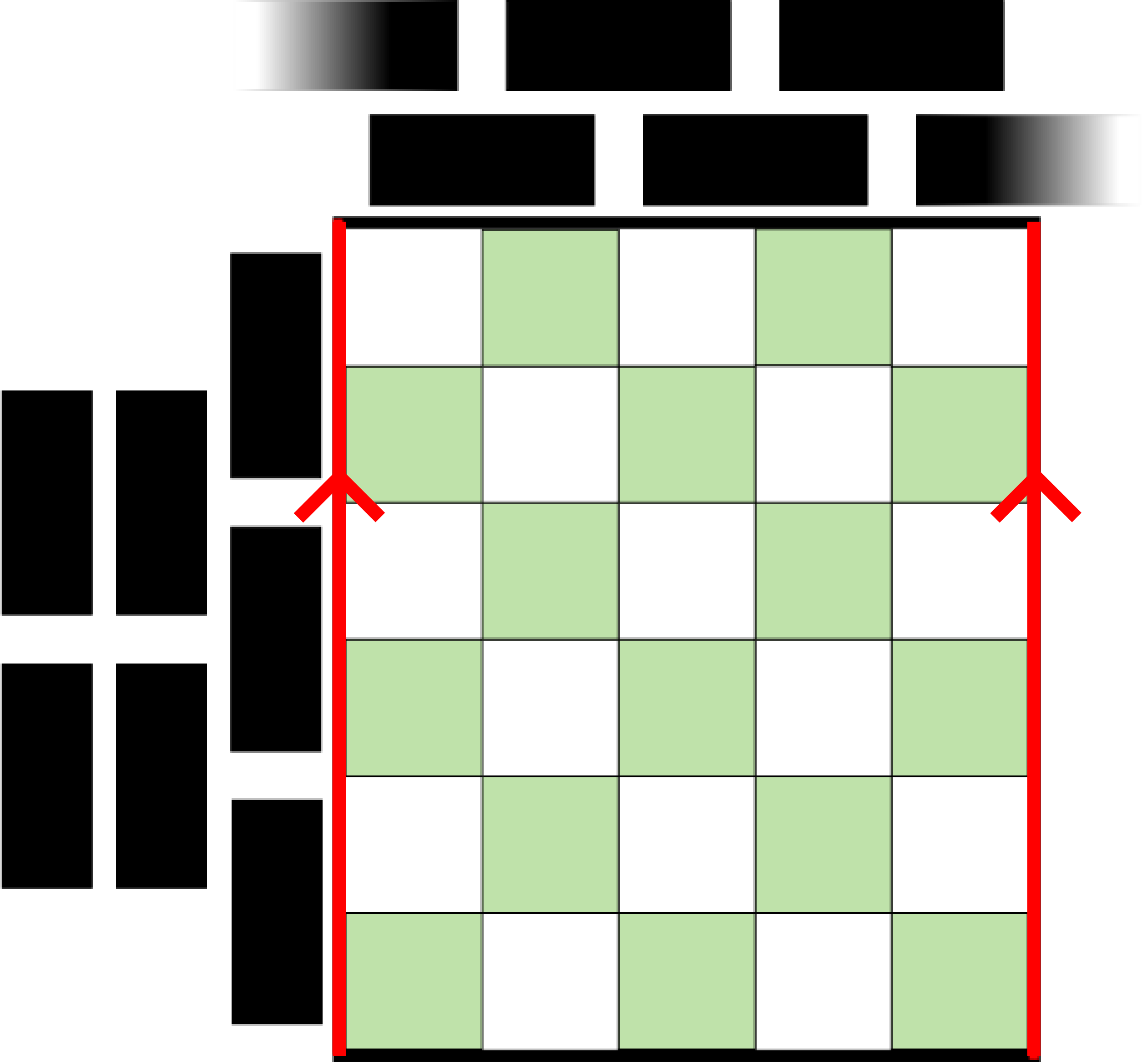}
         \caption{The $6^\prime \times 5$ board and its required tiles to be fault-free tileable}
         \label{fig:5x6'}
\end{figure}

Similarly, the $6^\prime \times 5$ board is not fault-free tileable. As shown in Figure~\ref{fig:5x6'}, this board requires $12$ tiles to cover all of its fault lines, but $15$ tiles are needed to fill the board. Again, tiles can be added in pairs or a new set of $5$. The additional $3$ tiles needed to fill the board cannot be added using sets of $2$ and $5$. Thus, the remaining $3$ tiles cannot be fit into the board, and it is not fault-free tileable.

\begin{figure}[h]
     \centering
         \includegraphics[scale=.25]{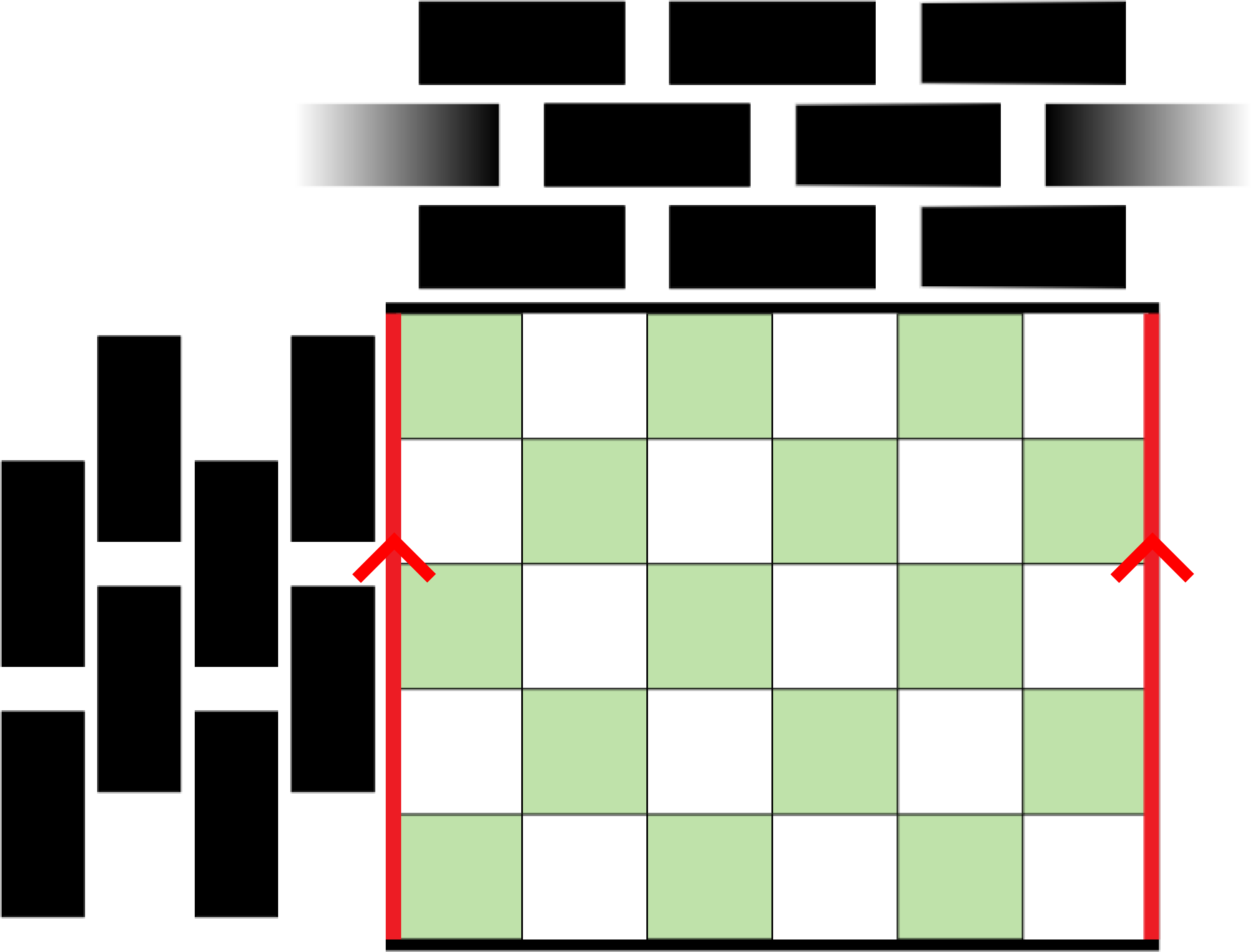}
         \caption{The $5^\prime \times 6$ board and its required tiles to be fault-free tileable}
         \label{fig:6x5'}
\end{figure}

The $5^\prime \times 6$ is not fault-free tileable. This board requires $17$ tiles to cover all of the fault-lines, as shown in Figure~\ref{fig:6x5'}. However, there is only room for $15$ .

\begin{figure}[h]
     \centering
         \includegraphics[scale=1]{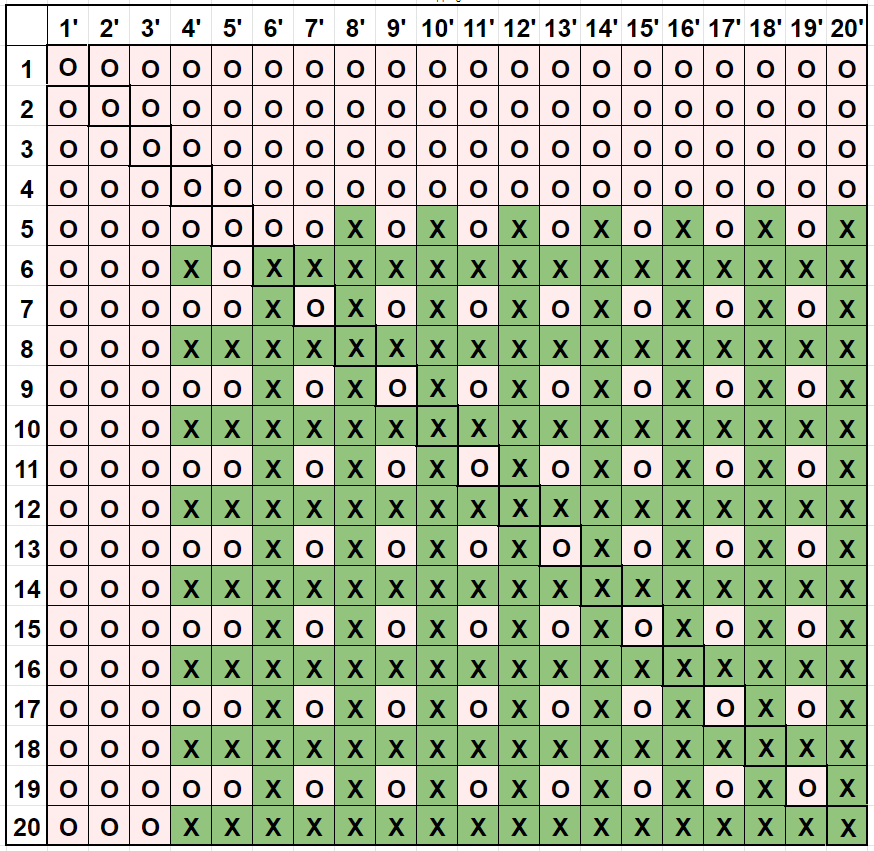}
         \caption{A chart denoting fault-free tileablitiy of cylindrical boards where boxes marked by an X are fault-free tileable, and boxes marked by an O are not}
         \label{fig:cylinderchart}
\end{figure}

The cases given above determine the fault-free tileability of all possible cylindrical boards. In Figure~\ref{fig:cylinderchart}, the fault-free tileability of cylinders up to $20^\prime \times 20$ are shown. Boxes marked by an X are fault-free tileable, while boxes marked by an O are not.

\section{Tori}

\begin{figure}[h]
     \centering
         \includegraphics[scale=.5]{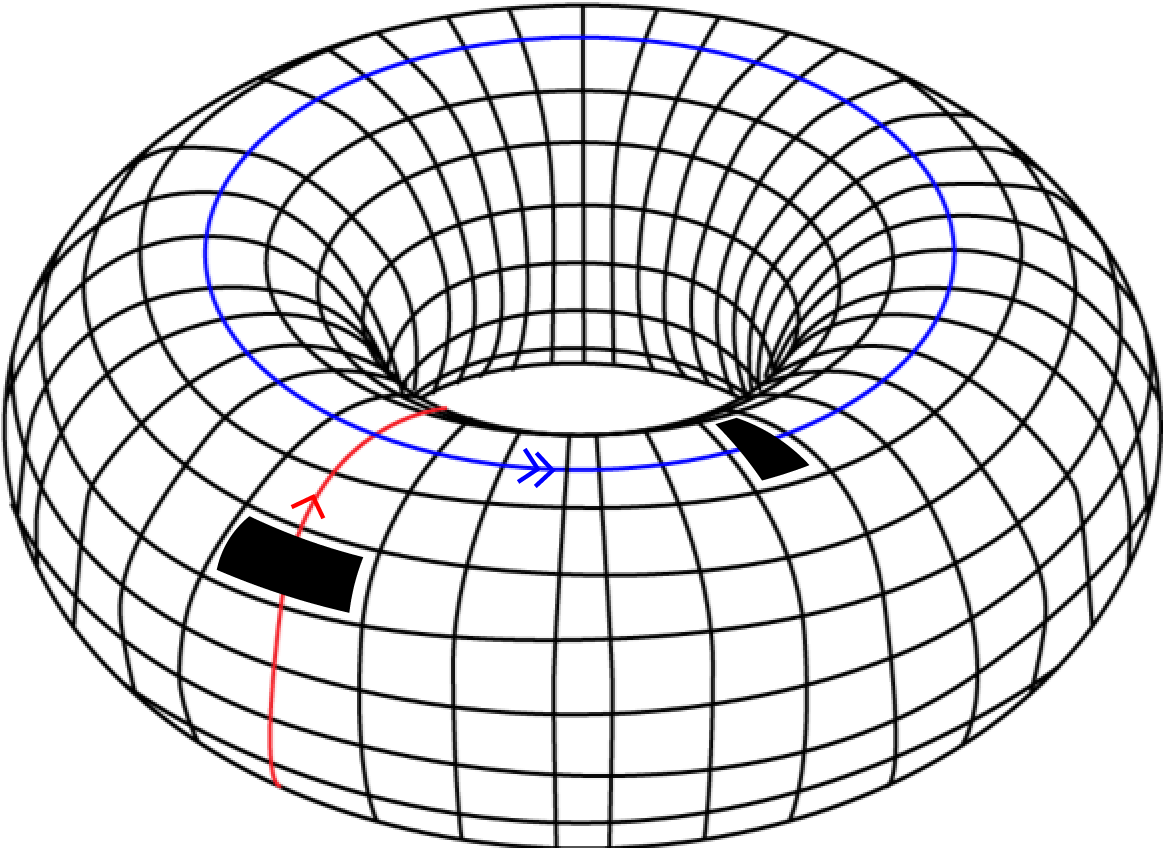}
         \caption{A $37^\prime \times 16^\prime$ torus board shown in three dimensions}
         \label{fig:literaltori}
\end{figure}

\begin{figure}[h]
     \centering
         \includegraphics[scale=.25]{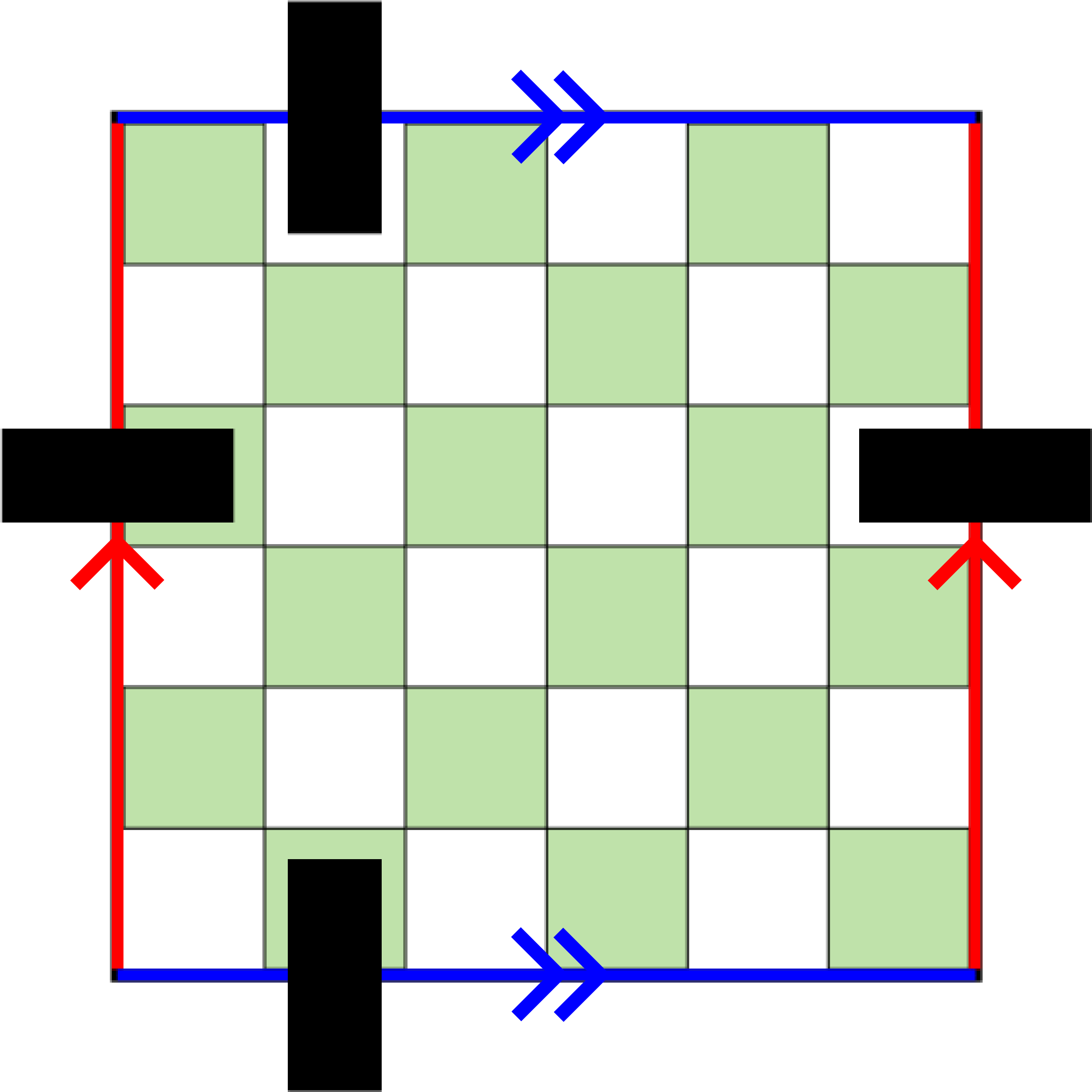}
         \caption{The $6^\prime \times6^\prime$ board with two tiles placed}
         \label{fig:torusex}
\end{figure}

Next, we turn to boards that are tori. A torus can be represented as a rectangular board where both sets of opposing edges are glued together. This is shown in Figure~\ref{fig:literaltori}. A torus board will have two more fault-lines than a rectangular board of the same dimensions. This is because both sets of edges in the rectangle have become fault-lines in the torus. A torus board is denoted as $a^\prime \times b^\prime $, where $a$ is the vertical length and $b$ the horizontal length, as shown in Figure~\ref{fig:torusex}. Both $a$ and $b$ are denoted with primes because in both directions the board wraps around and connects to itself.

\begin{figure}[h]
\centering
         \subcaptionbox{\label{fig:4'x4'}}{\includegraphics[scale=.25]{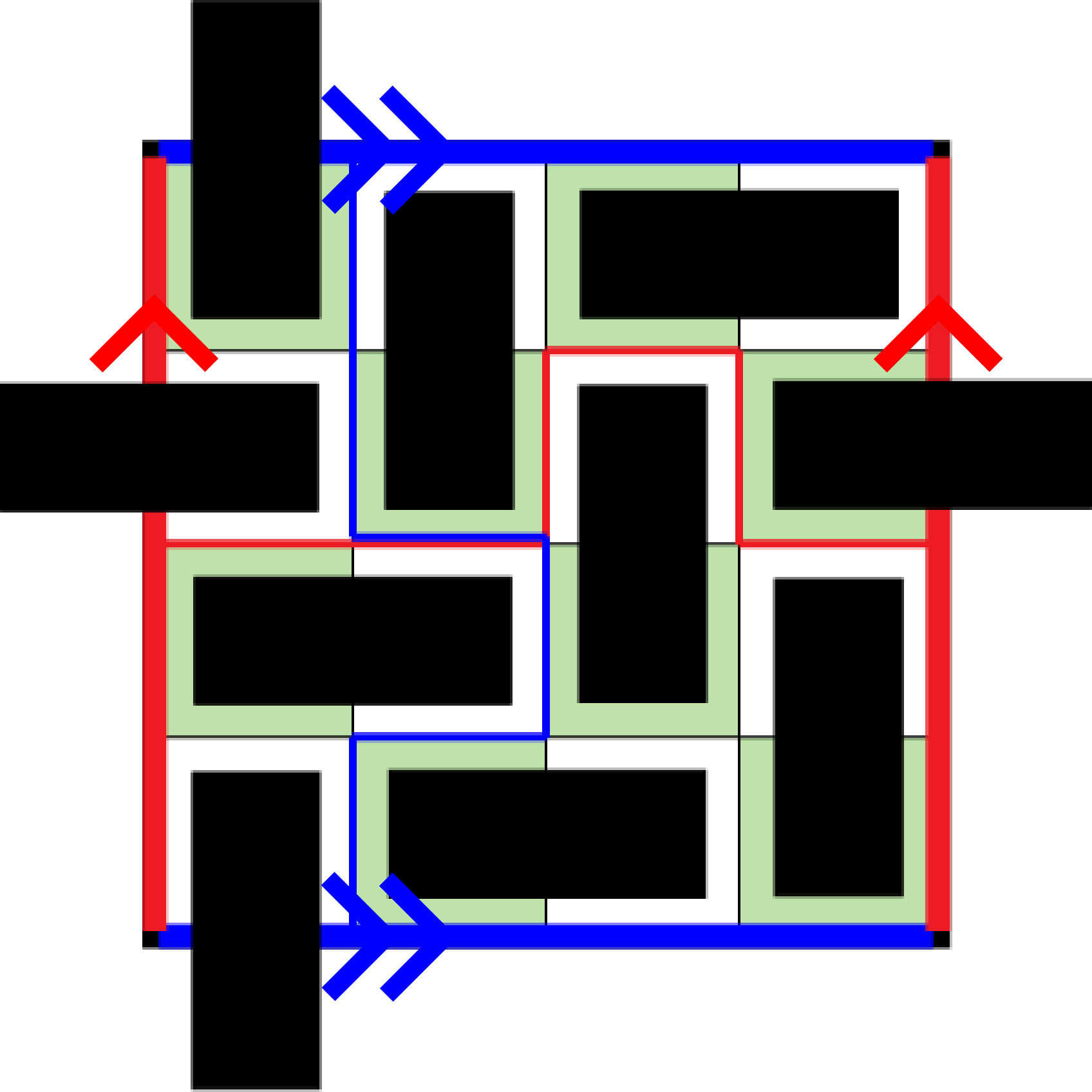}}
         \vspace*{.5cm}
         \subcaptionbox{\label{fig:7'x8'}}{\includegraphics[scale=.25]{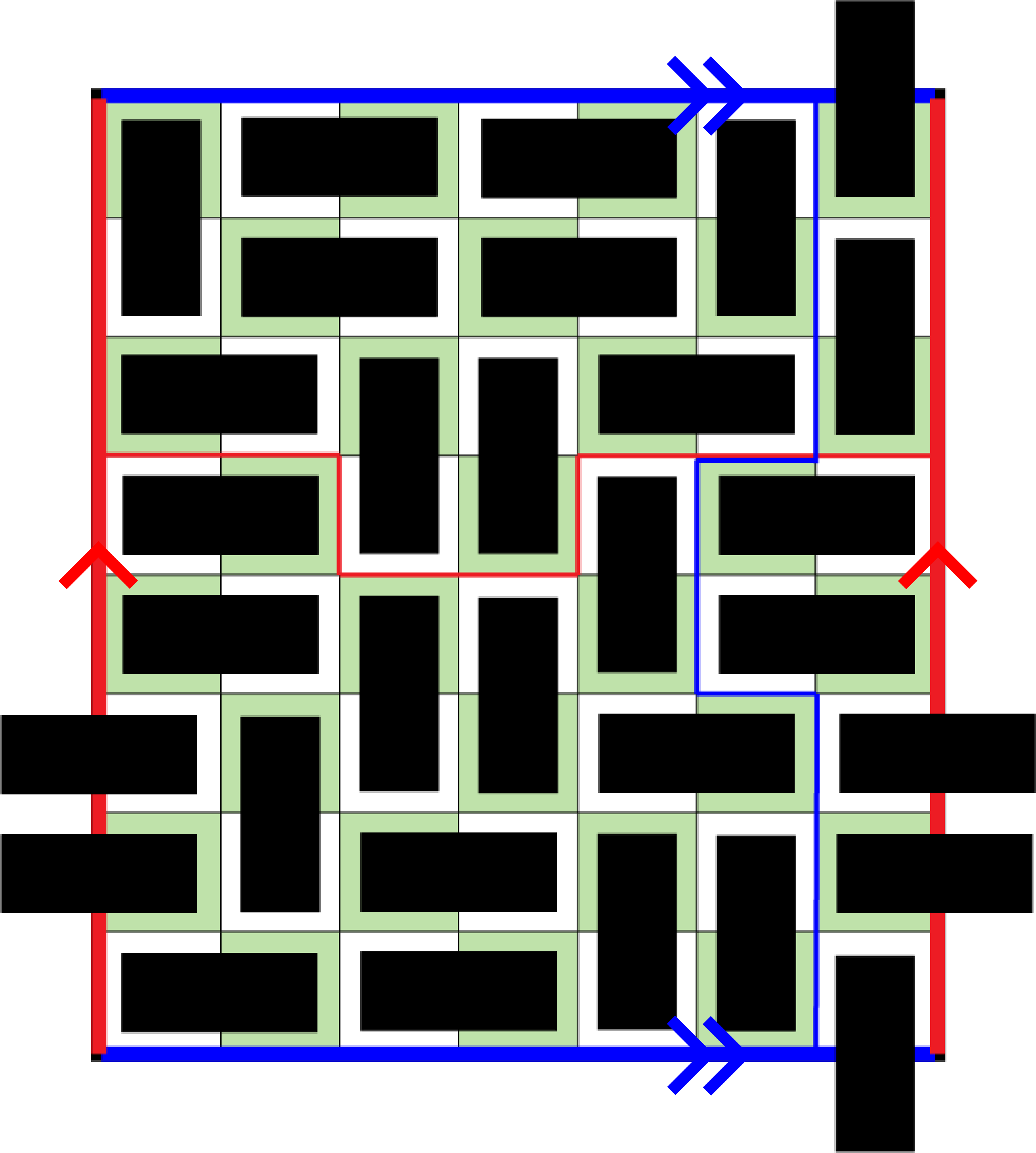}}
         
         \hspace*{1cm}
          \subcaptionbox{\label{fig:6'x9'}}{\includegraphics[scale=.25]{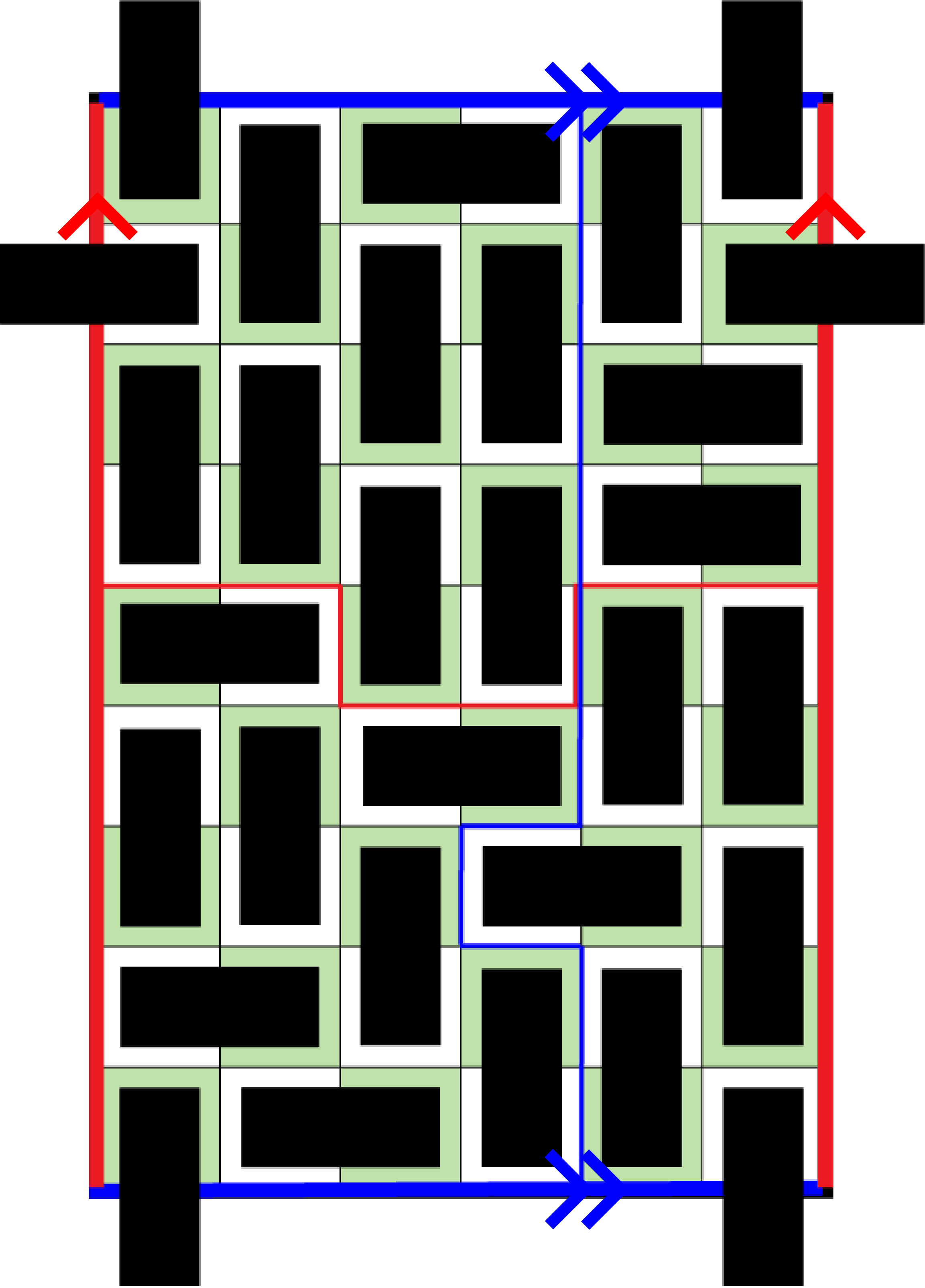}}
          \vspace*{.5cm}
          \subcaptionbox{\label{fig:5'x10'}}{\includegraphics[scale=.25]{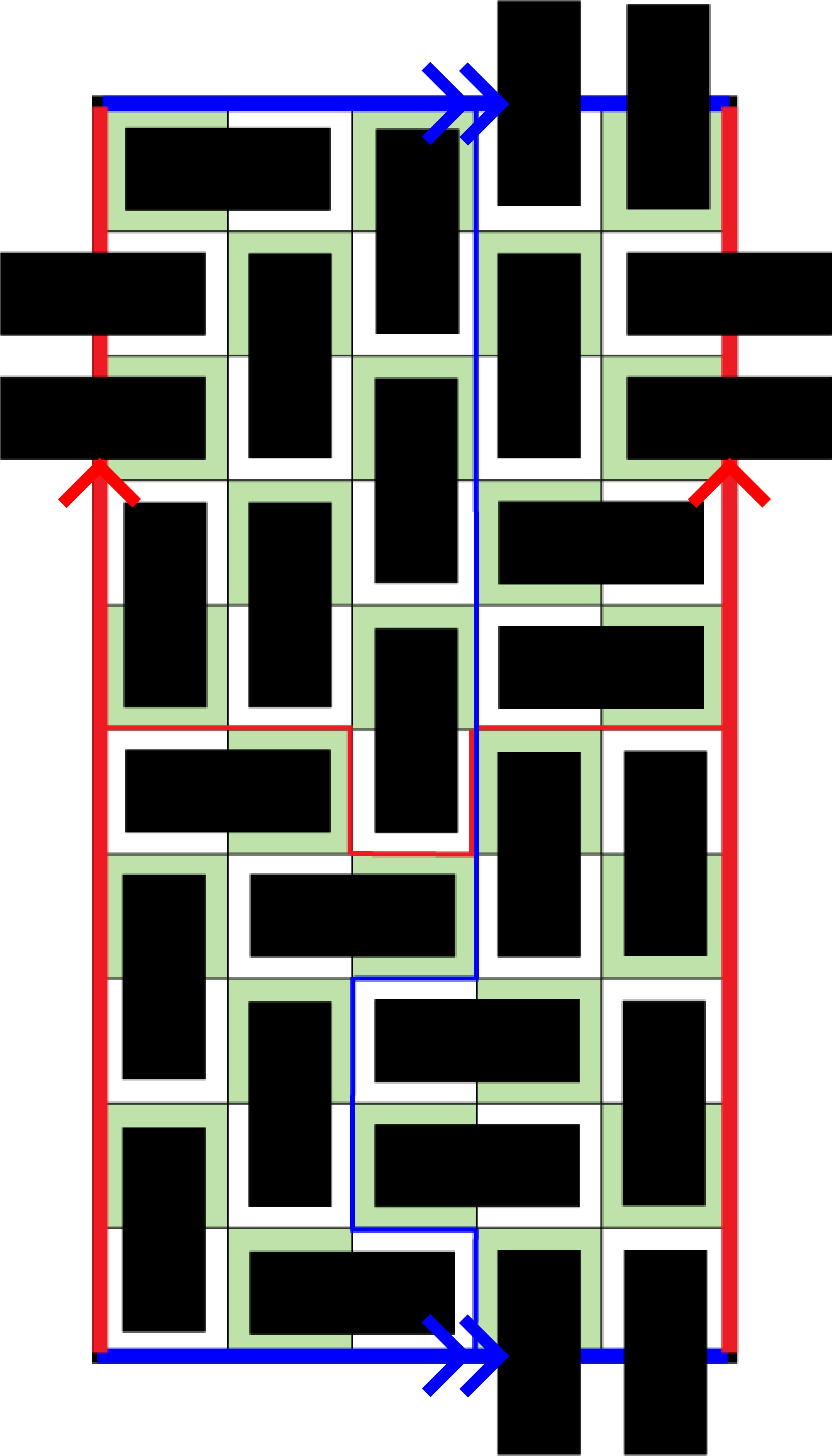}}\vspace*{.5cm}
        
        \caption{The fault-free tileable torus boards}
        \label{fig:tileabletori}
\end{figure}

The following tori are fault-free tileable as shown in Figure~\ref{fig:tileabletori}:
\begin{enumerate}[(a)]
    \item $4^\prime \times 4^\prime$
    \item $8^\prime \times 7^\prime$
    \item $9^\prime \times 6^\prime$
    \item $10^\prime \times 5^\prime$
\end{enumerate}
Again, these can be expanded by multiples of two in either direction. Thus the following are fault-free tileable, where $n,m \in \mathbb{N}$:
\begin{enumerate}[(a)]
    \item $(4+2n)^\prime \times (4+2m)^\prime$
    \item $(8+2n)^\prime \times (7+2m)^\prime$
    \item $(9+2n)^\prime \times (6+2m)^\prime$
    \item $(10+2n)^\prime \times (5+2m)^\prime$
\end{enumerate}

We now show why the following boards are not fault-free tileable.
\begin{itemize}
    \item $\text{odd}^\prime \times \text{odd}^\prime$
    \item $(2n)^\prime \times 1^\prime $
    \item $(2+n)^\prime \times 2^\prime $
    \item $(4+2n)^\prime \times 3^\prime $
    \item $(5+2n)^\prime \times 4^\prime $
    \item $6^\prime \times 5^\prime $
    \item $8^\prime \times 5^\prime $
    \item $7^\prime \times 6^\prime $
\end{itemize}

The odd by odd torus boards are not fault-free tileable, as before, by parity.

The $(2n)^\prime \times 1^\prime $ boards are not fault-free tileable. Any first tile placement will create a fault-line.

The $(2+n)^\prime \times 2^\prime $ boards are not fault-free tileable. A vertical fault-line must be crossed, thus a horizontal tile must be placed to cover it. This will result in a horizontal fault-line above and below the tile.

Similarly, the $(4+2n)^\prime \times 3^\prime $ boards are not fault-free tileable. A vertical fault-line must be crossed, thus a horizontal tile must be placed to cover it. A vertical tile must then be placed in the single remaining space of the row, which will result in a horizontal fault-line.

\begin{figure}[h]
     \centering
         \includegraphics[scale=.25]{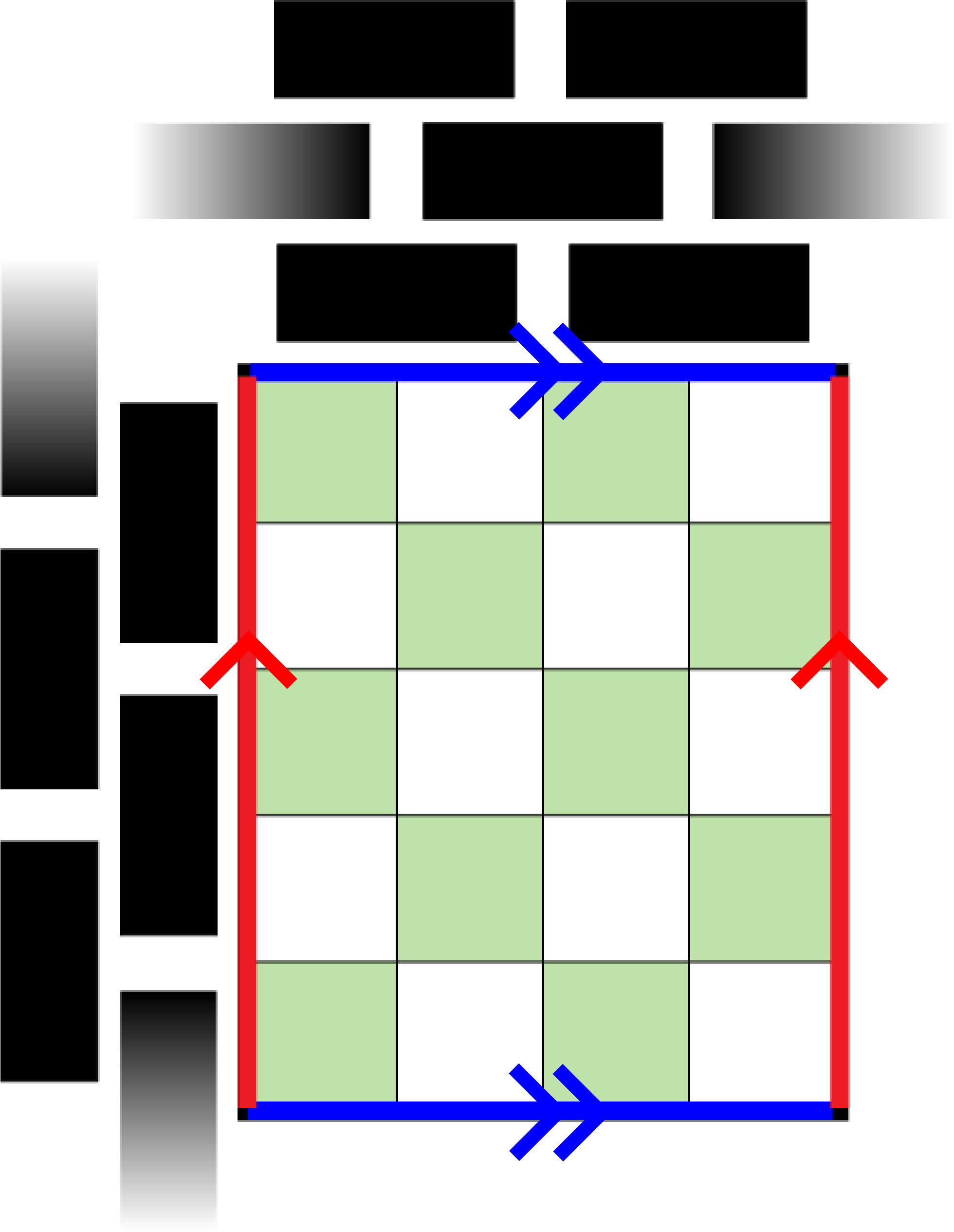}
         \caption{The $5^\prime \times 4^\prime$ and its required tiles to be fault-free tileable}
         \label{fig:4'x5'}
\end{figure}

The $(5+2n)^\prime \times 4^\prime $ boards are not fault-free tileable. An example on the $5^\prime \times 4^\prime$ board is shown in Figure~\ref{fig:4'x5'}. These boards require 6 tiles to cover all vertical fault-lines and $5+2n$ tiles to cover the horizontal fault-lines. This is a total of $11+2n$ tiles, an odd number. A $(5+2n)^\prime \times 4^\prime $ board tiled fully will have $10+4n$ tiles, which is an even number of tiles. We can add tiles to the board in pairs, a new set of $6$, or a new set of $5+2n$. Similarly to the $4^\prime \times (5+2n)$ boards, any other set of tiles will create a row or column with an odd number of remaining spaces to be covered. To get an even number of tiles to fill the board, we must add a new set of $5+2n$. This will give us $16+4n$ tiles, which is more than there is space for on the board.

\begin{figure}[h]
     \centering
         \includegraphics[scale=.23]{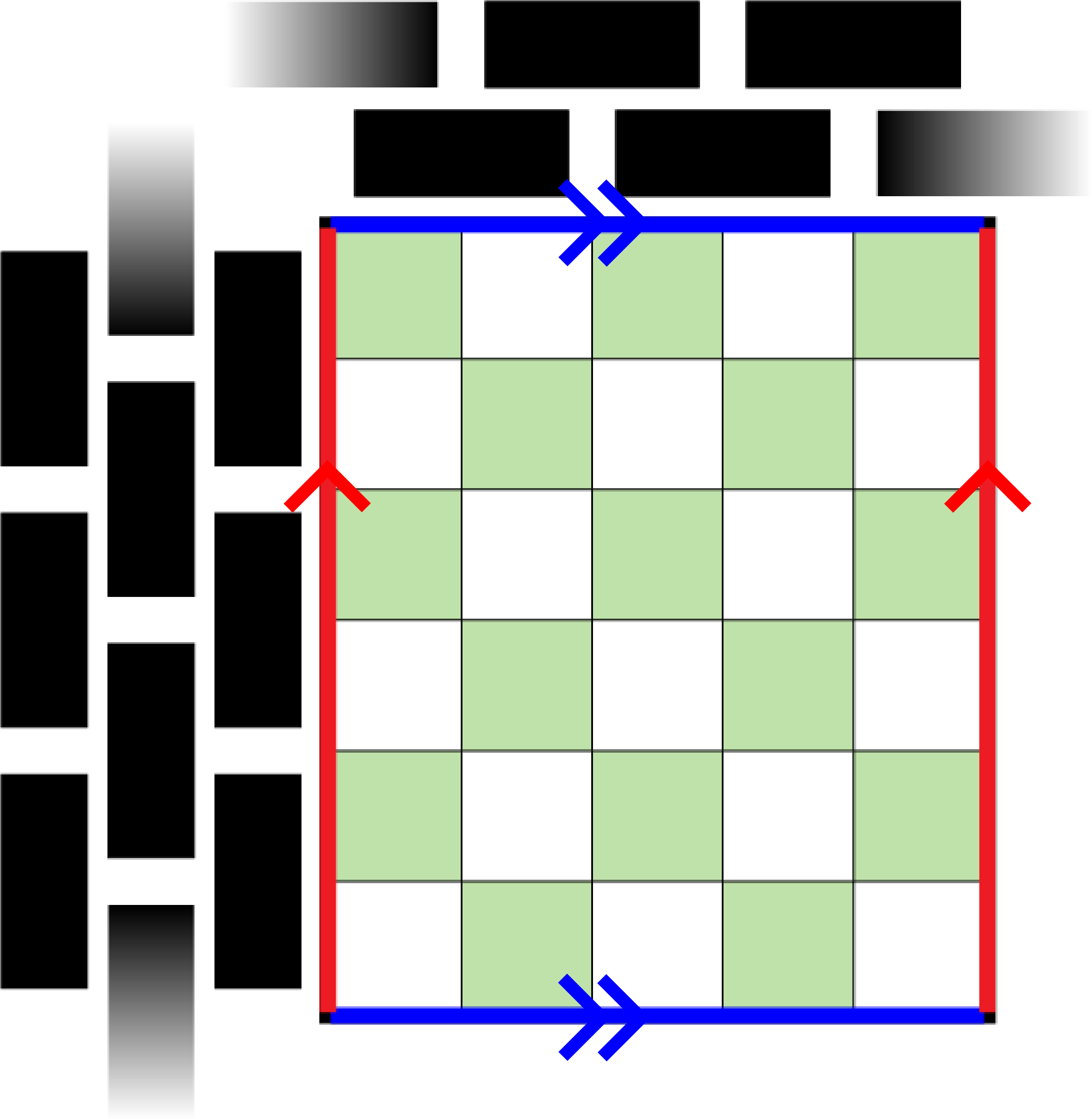}
         \caption{The $6^\prime \times 5^\prime$ and its required tiles to be fault-free tileable}
         \label{fig:5'x6'}
\end{figure}

The $6^\prime \times 5^\prime $ board is not fault-free tileable. This requires $14$ tiles to cover all of the fault-lines, as shown in Figure~\ref{fig:5'x6'}. There is room for $15$ tiles. The last required tile cannot be added to the board, as tiles can only be added in pairs or a new set of $5$.

\begin{figure}[h]
\centering
         \subcaptionbox{\label{fig:5'x8'}}{\includegraphics[scale=.25]{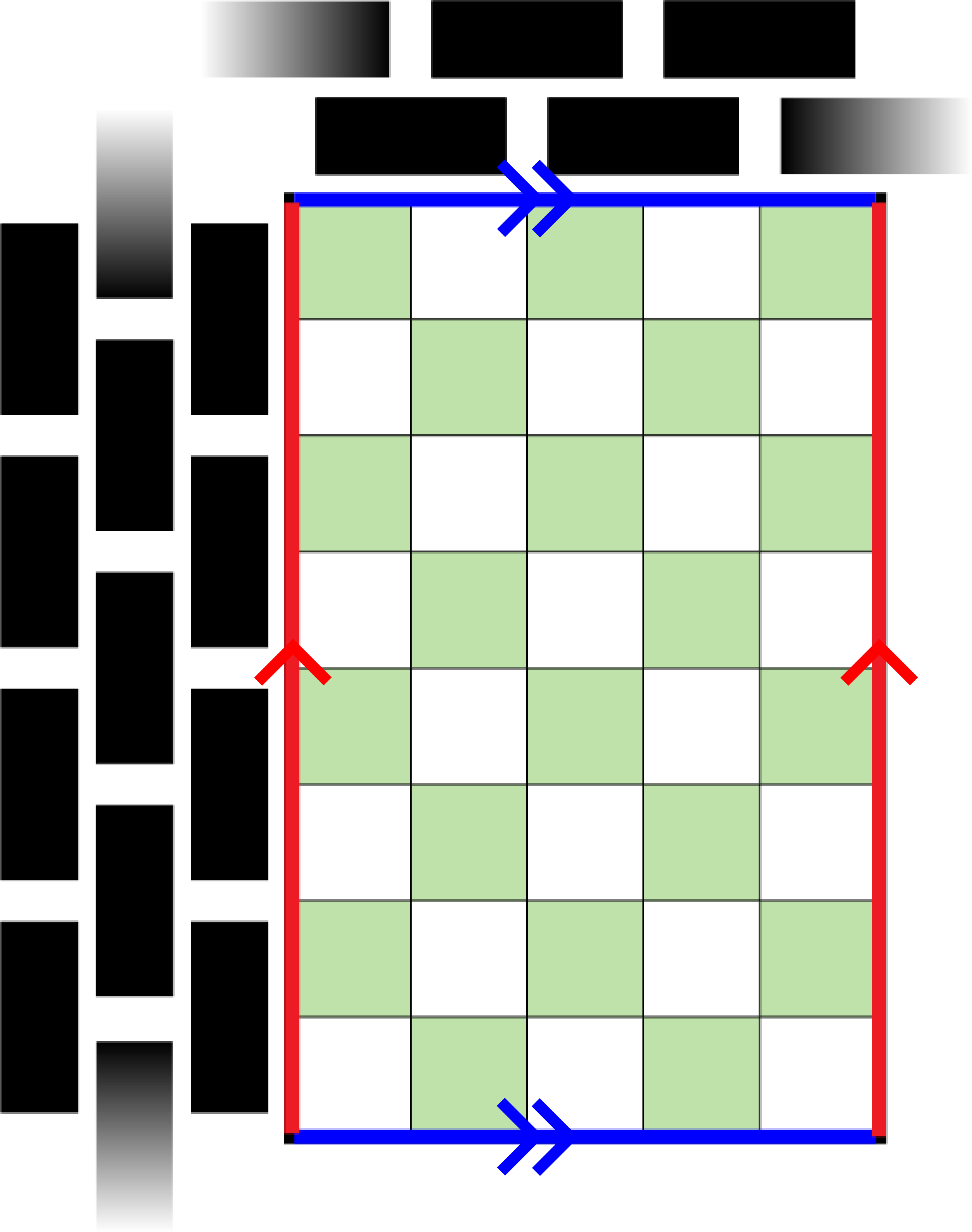}}
         \hspace*{1cm}
         \subcaptionbox{\label{fig:6'x7'}}{\includegraphics[scale=.25]{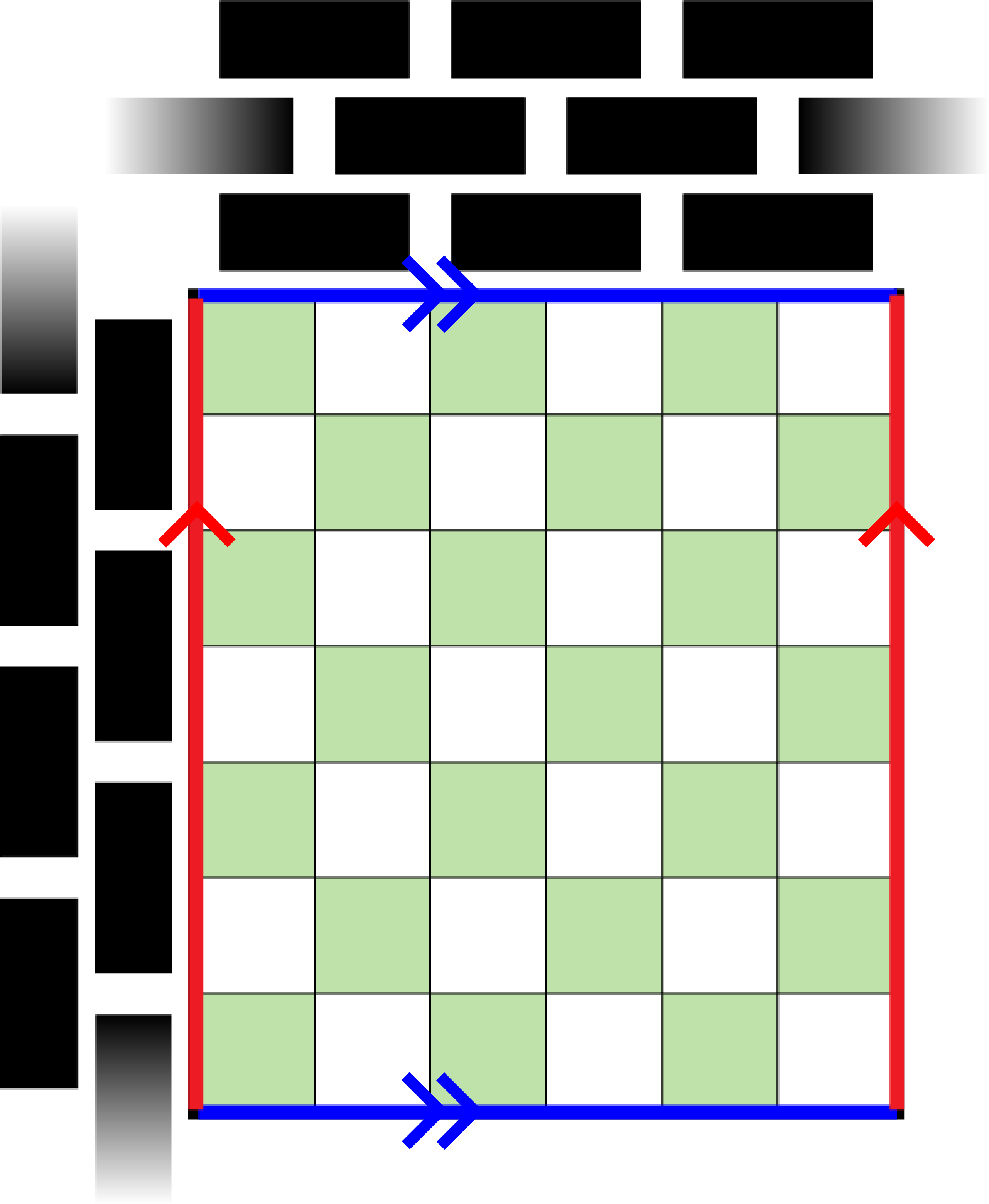}}
       
        \caption{The $8^\prime \times 5^\prime$ and $7^\prime \times 6^\prime$ and their required tiles to be fault-free tileable}
        \label{fig:5'x8'and6'x7'}
\end{figure}

Similarly, the $8^\prime \times 5^\prime $ and $7^\prime \times 6^\prime $ boards are not fault-free tileable. Their required tiles are shown in Figure~\ref{fig:5'x8'and6'x7'}. Again more tiles are required than the number of spaces allows for.

\begin{figure}[h]
     \centering
         \includegraphics[scale=1]{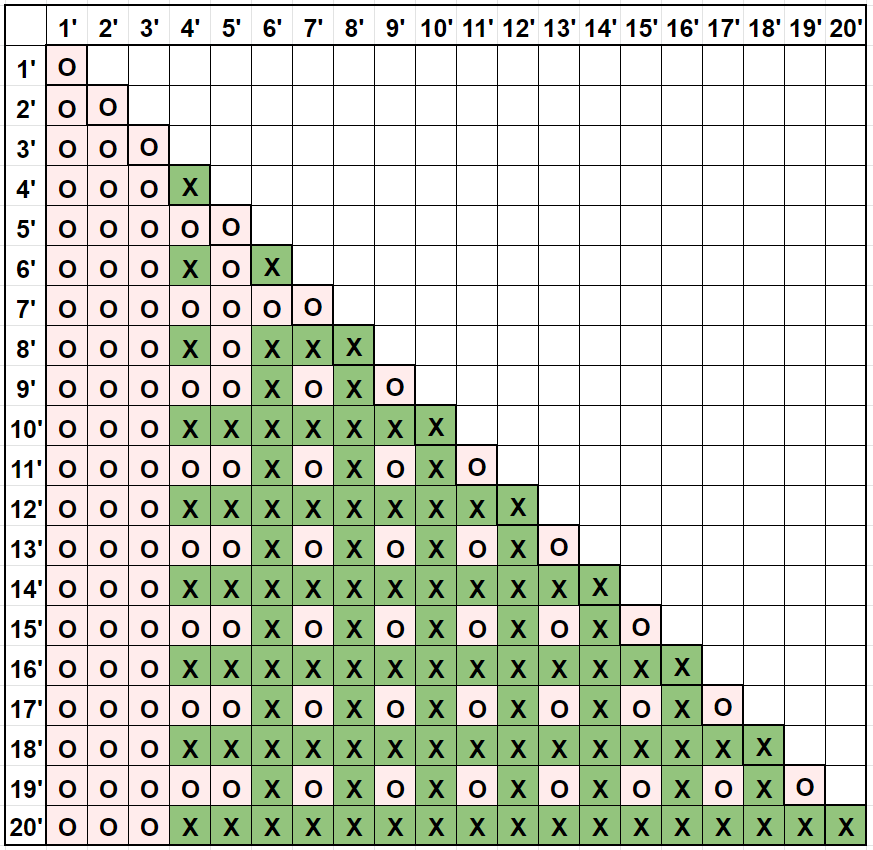}
         \caption{A chart denoting fault-free tileablitiy of torus boards where boxes marked by an X are fault-free tileable, and boxes marked by an O are not}
         \label{fig:toruschart}
\end{figure}

The cases above determine the fault-free tileabilty of all torus boards. In Figure~\ref{fig:toruschart} tileability of tori up to $20 ^\prime \times 20 ^\prime$ are shown.

\section{M\"{o}bius Strips}

\begin{figure}[h]
     \centering
         \includegraphics[scale=.5]{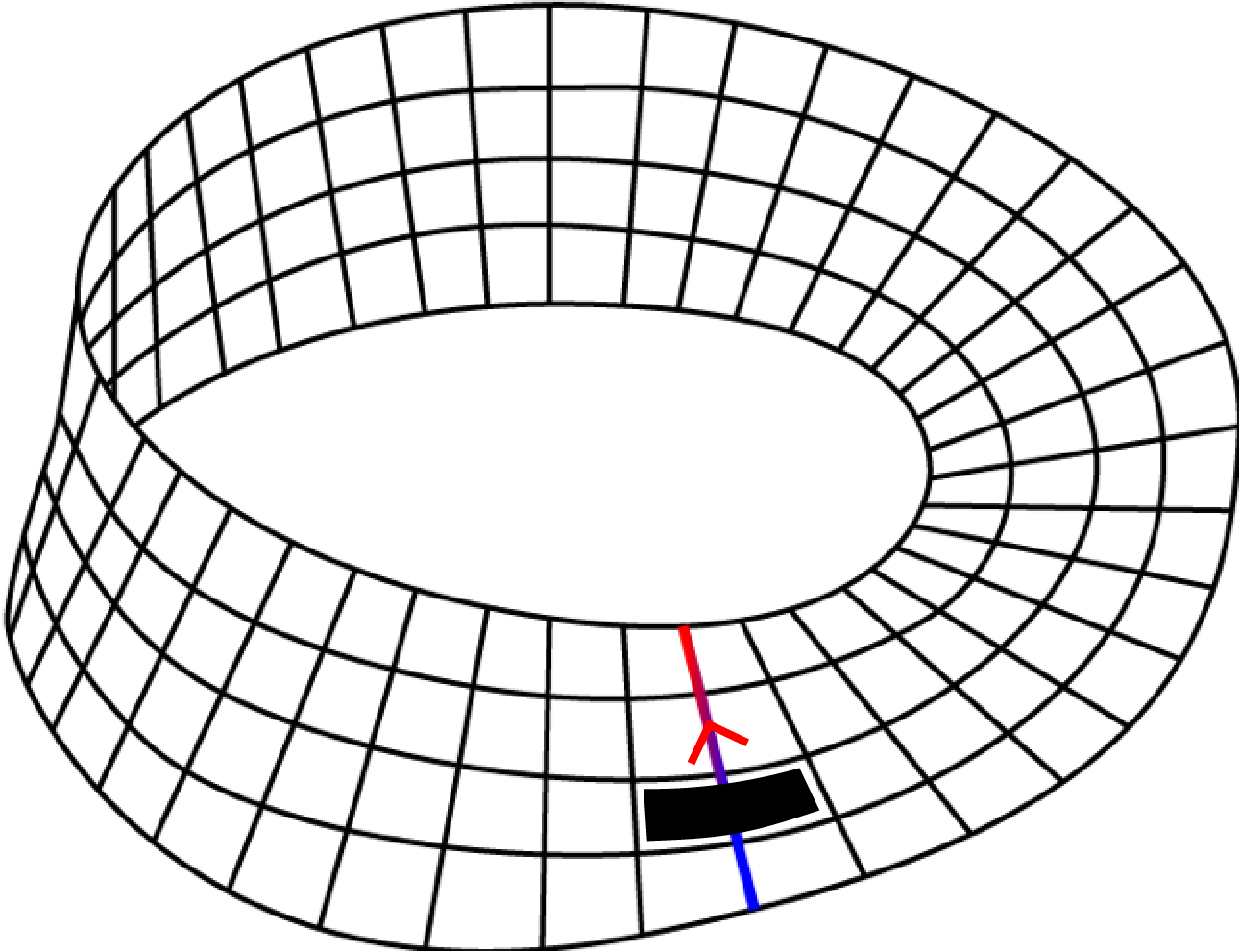}
         \caption{A $4^{\prime \prime} \times 38$ M\"{o}bius strip board shown in three dimensions}
         \label{fig:literalmobius}
\end{figure}

\begin{figure}[h]
\centering
         \subcaptionbox{\label{fig:mobiusex1}}{\includegraphics[scale=.25]{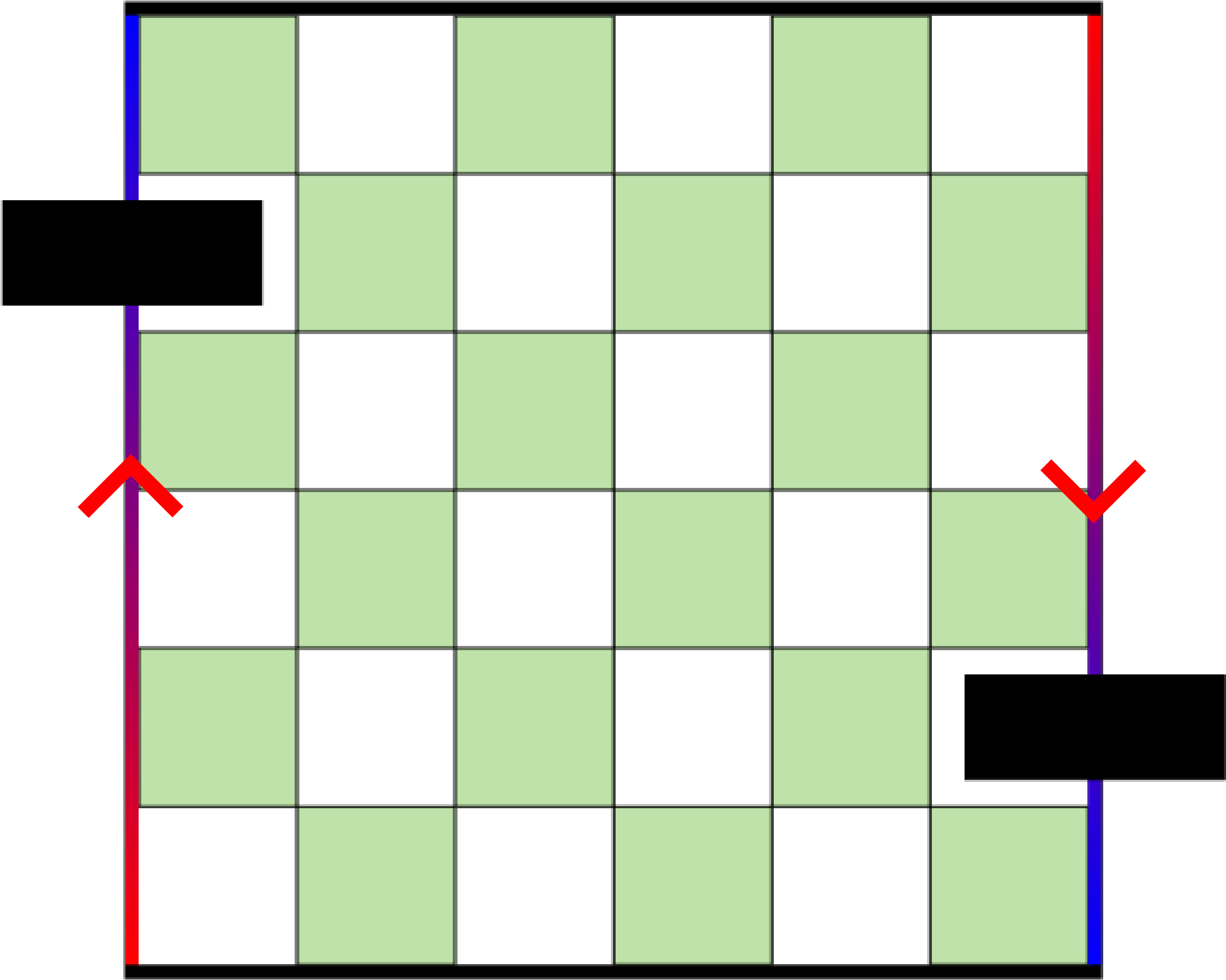}}
         \hspace*{1cm}
         \subcaptionbox{\label{fig:mobiusex2}}{\includegraphics[scale=.25]{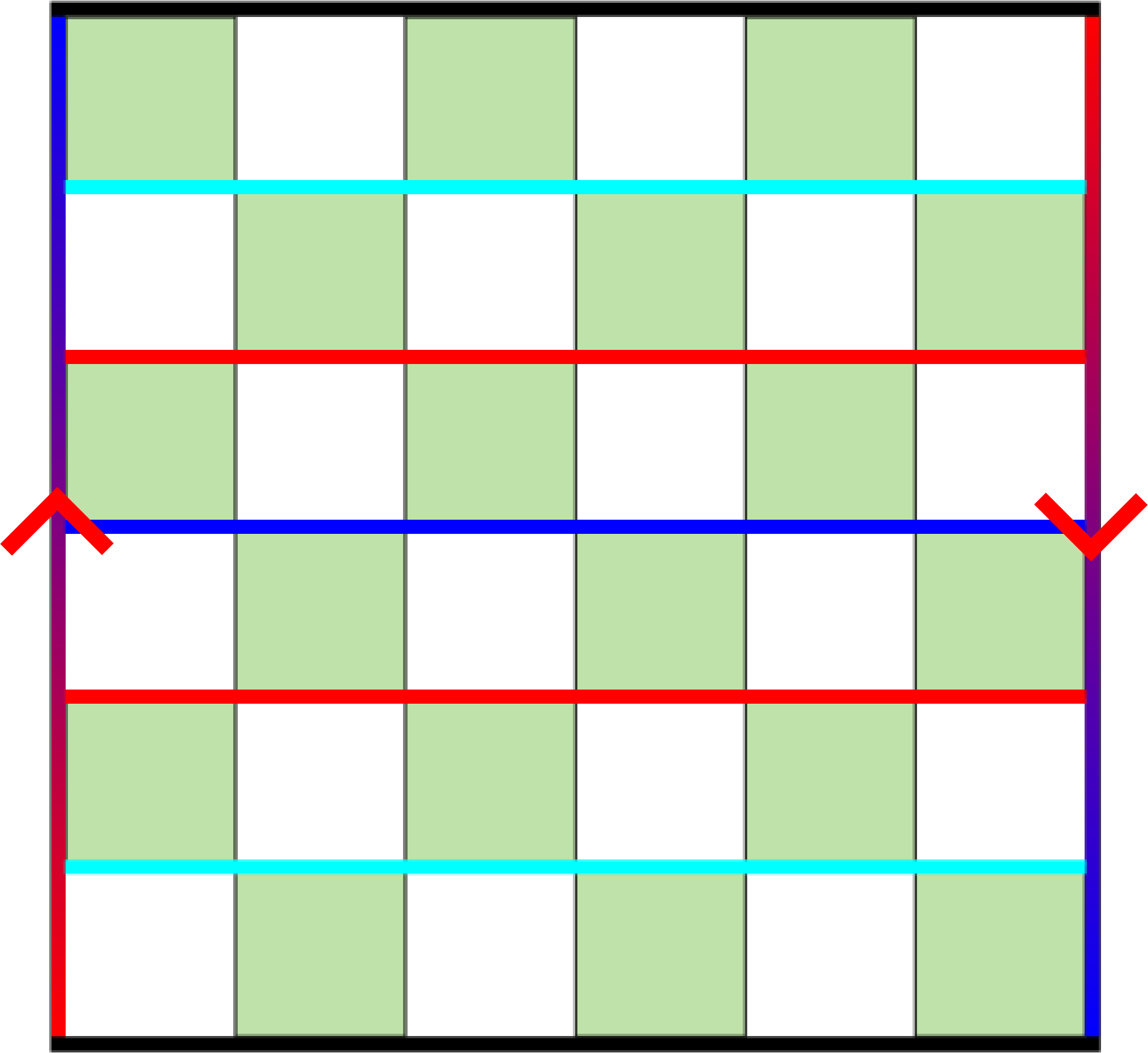}}
       
        \caption{A M\"{o}bius board with one tile placed and distinguished fault-lines}
        \label{fig:mobiusexs}
\end{figure}

Lastly, we looked at the fault-free tileablity of M\"{o}bius strips. These are similar to cylinders. However, when imagining the edges being glued together, we introduce a half twist, as shown in Figure~\ref{fig:literalmobius}. These are denoted as $a^{\prime \prime} \times b$, where $a$ is the height or length of the vertical sides when shown in two dimensions, and $b$ is the distance around the M\"{o}bius strip or length of horizontal sides when shown in two dimensions. A side of length $a^{\prime \prime}$ represents a fault-line which glues the two opposing edges together after introducing a half-twist. One important difference in the tiling of M\"{o}bius strips is that the fault-lines also wrap around the twist, as shown in Figure~\ref{fig:mobiusexs}. Most fault-lines appear twice on the two dimensional image, though they only must be covered once. The exception to this is the center most horizontal fault-line on boards with an even height. These fault-lines appear only once and only must be covered once.

\begin{figure}[h]
\centering
         \subcaptionbox{\label{fig:3x4''}}{\includegraphics[scale=.25]{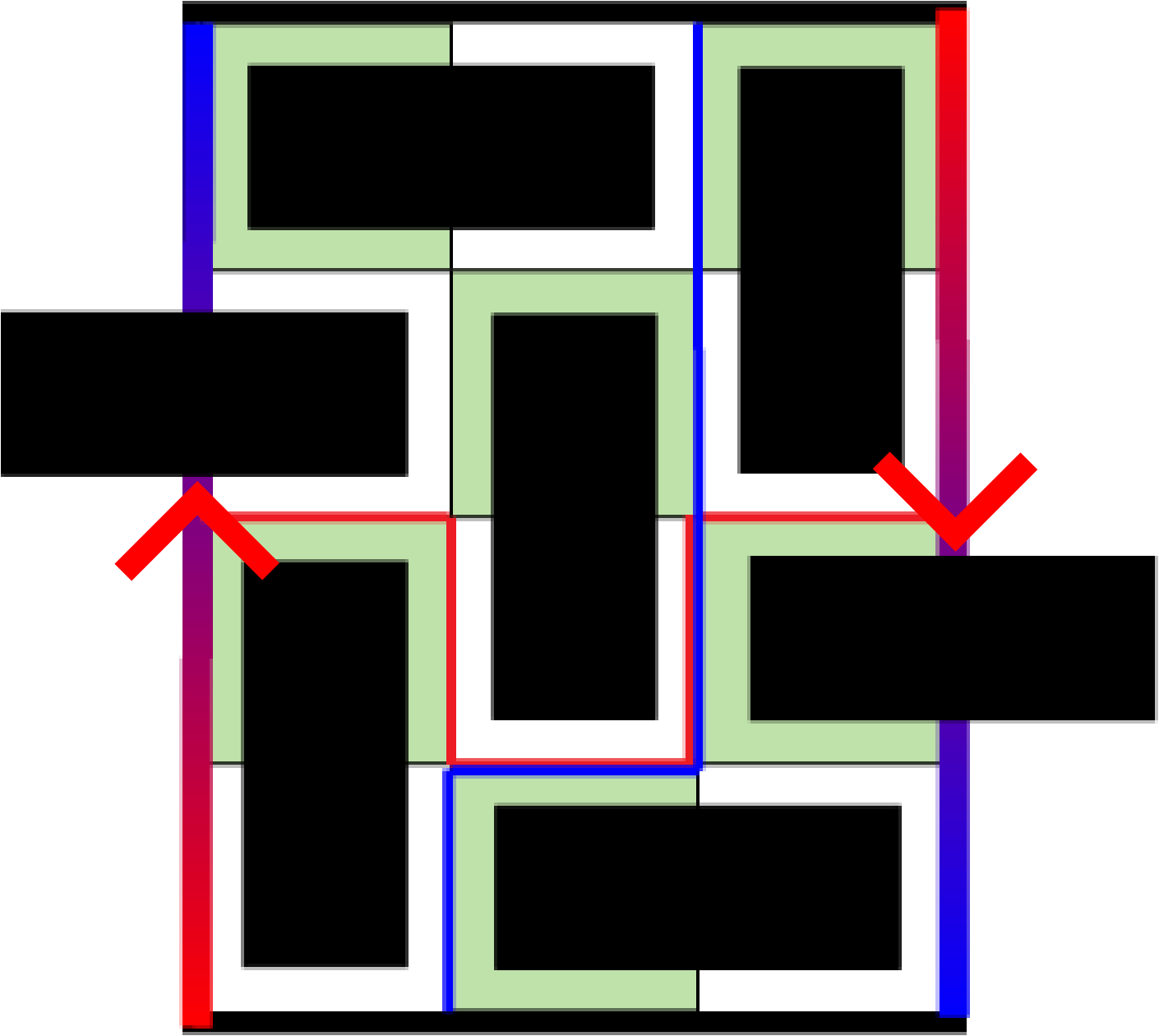}}
         
         \hspace*{1cm}
         \subcaptionbox{\label{fig:4x5''}}{\includegraphics[scale=.25]{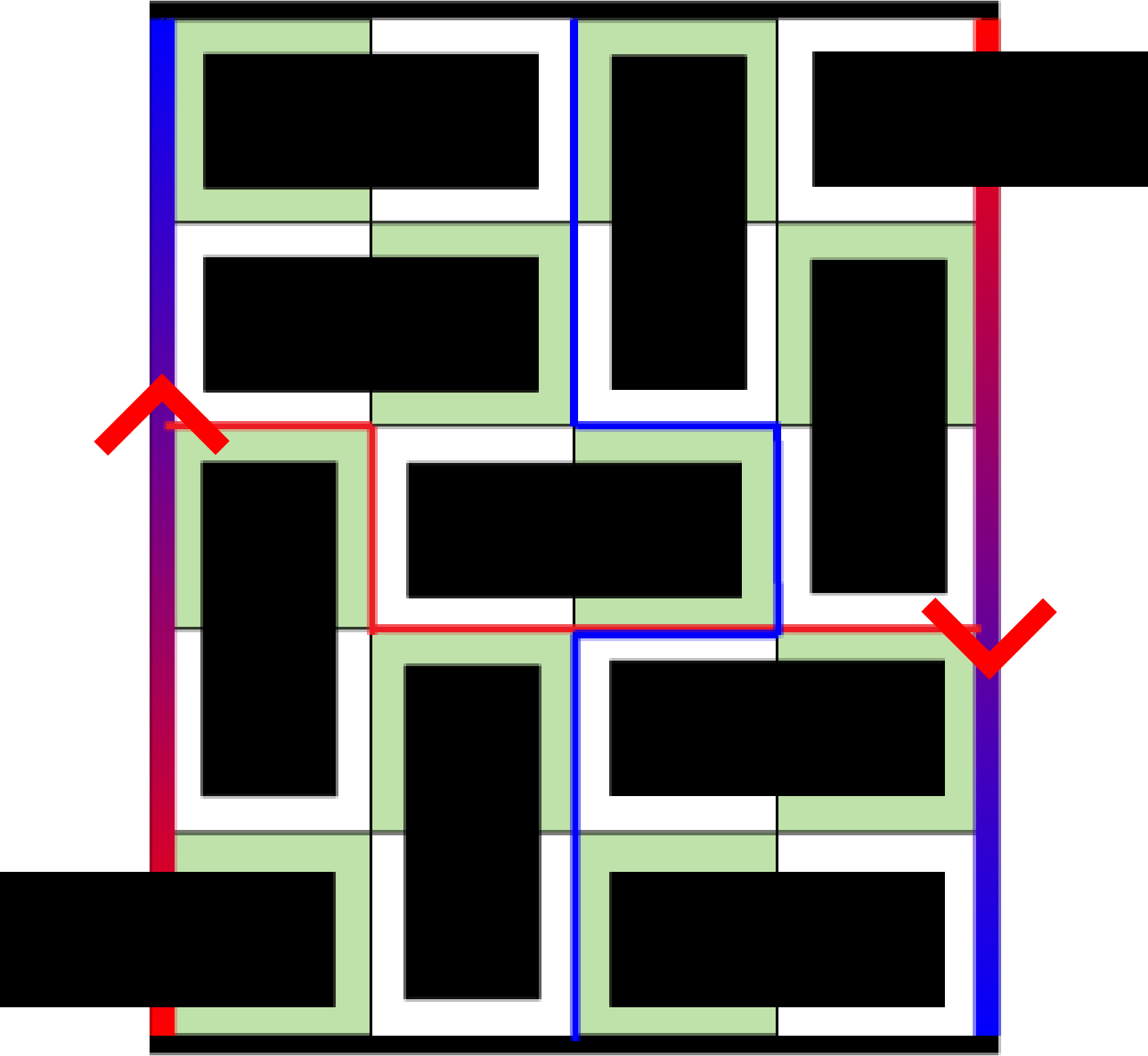}}\vspace*{.5cm}
          \subcaptionbox{\label{fig:5x4''}}{\includegraphics[scale=.25]{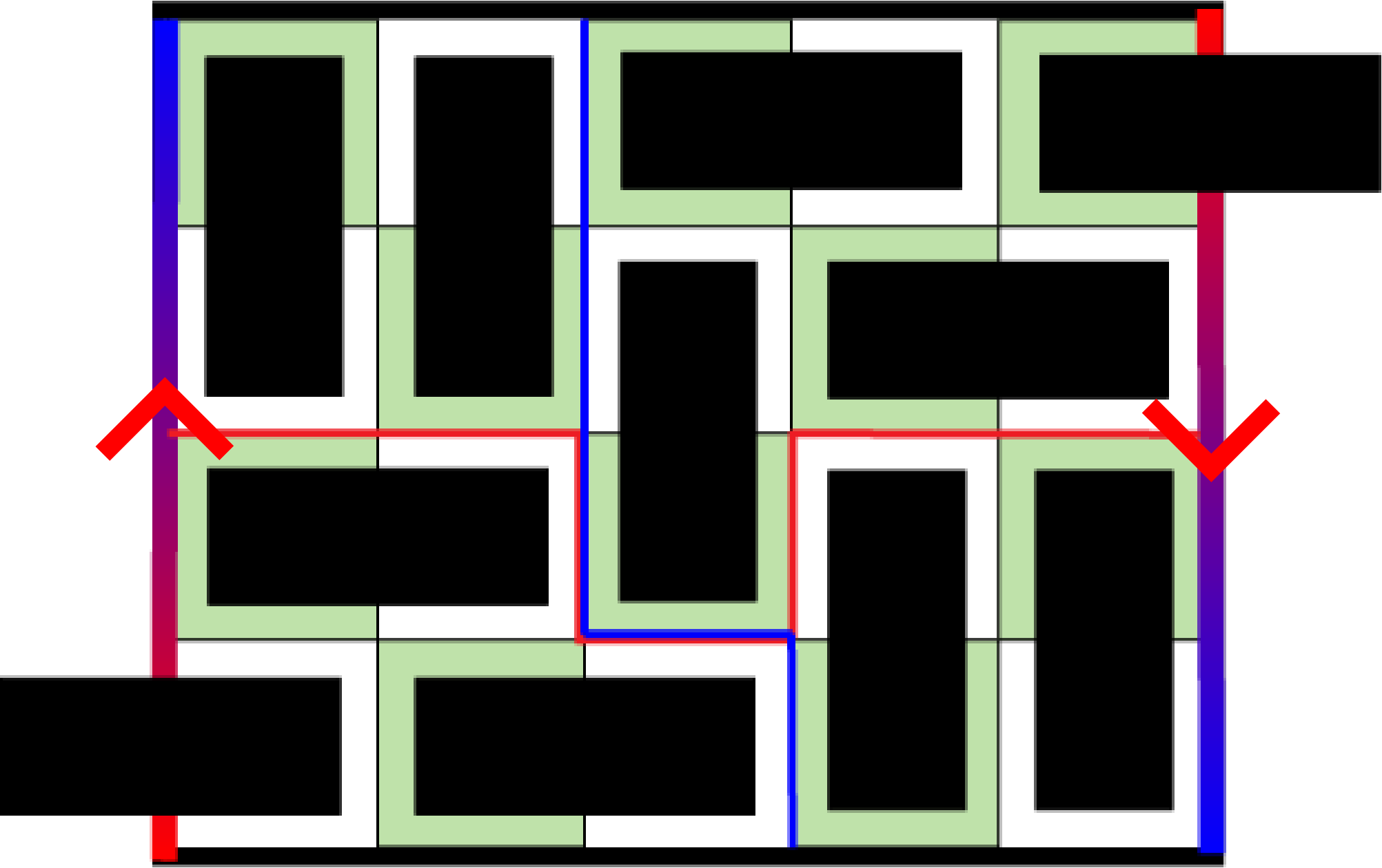}}
          
         \hspace*{1cm}
          \subcaptionbox{\label{fig:6x6''}}{\includegraphics[scale=.25]{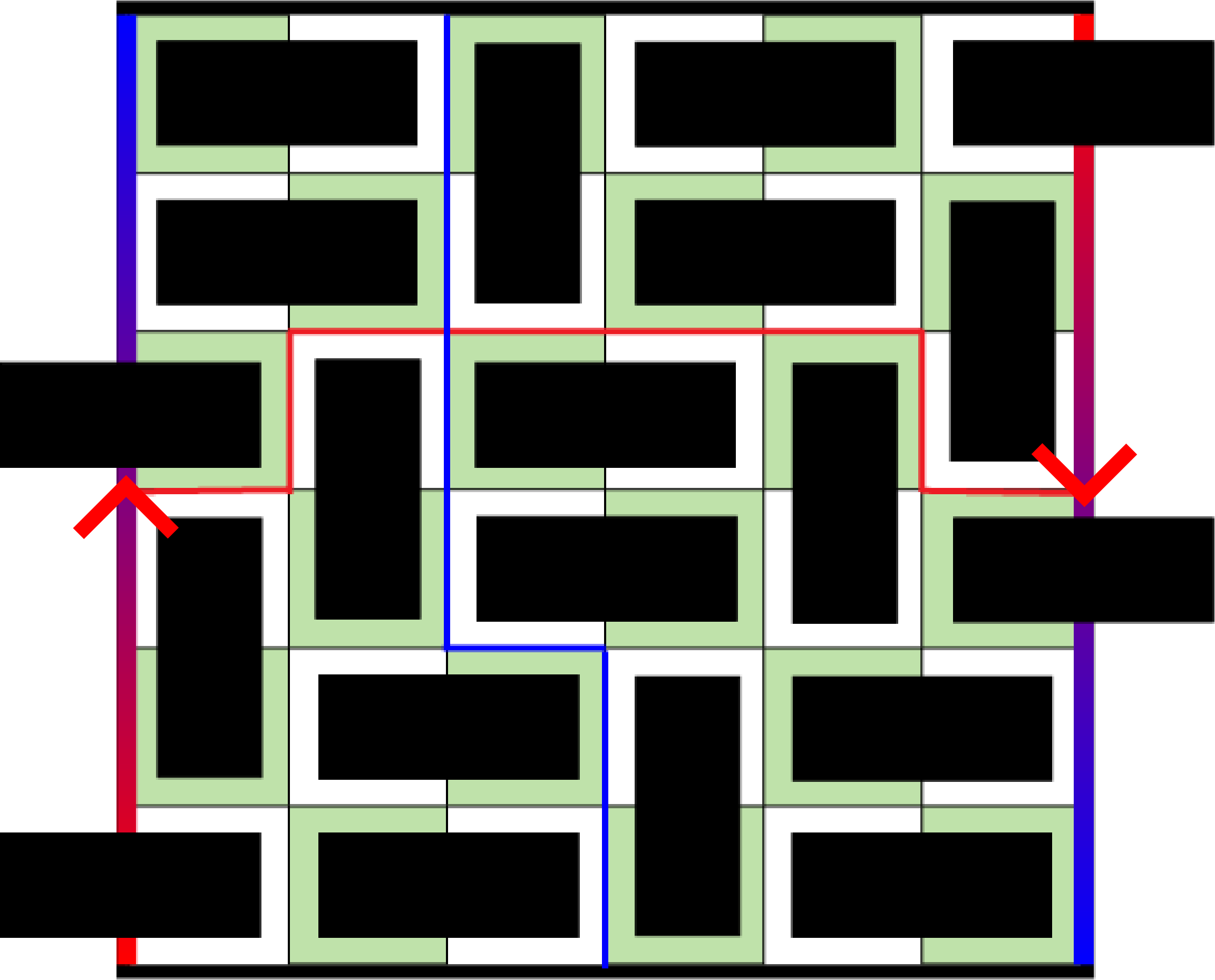}}\vspace*{.5cm}
         \subcaptionbox{\label{fig:4x8''}}{\includegraphics[scale=.25]{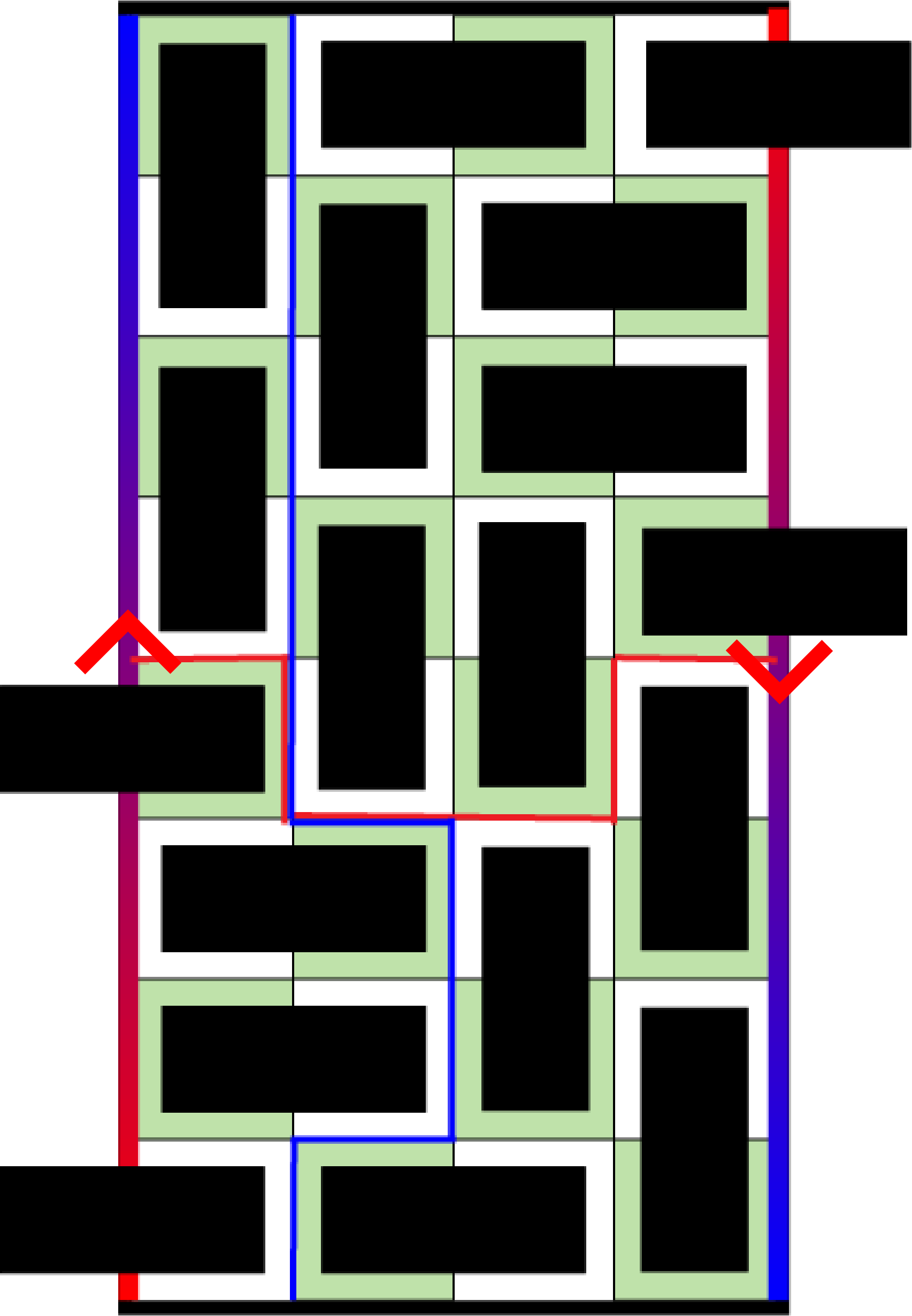}}

        \caption{The fault-free tileable M\"{o}bius boards}
        \label{fig:mobiusboards}
\end{figure}

The following M\"{o}bius strips are fault-free tileable as shown in Figure~\ref{fig:mobiusboards}:
\begin{enumerate}[(a)]
    \item $4^{\prime \prime} \times 3$
    \item $5^{\prime \prime} \times 4$
    \item $4^{\prime \prime} \times 5$
    \item $6^{\prime \prime} \times 6$
    \item $8^{\prime \prime} \times 4$
\end{enumerate}
Again, these can be expanded by multiples of two in either direction. Thus the following are fault-free tileable, where $n,m \in \mathbb{N}$:
\begin{enumerate}[(a)]
    \item $(4+2n)^{\prime \prime} \times (3+2m)$
    \item $(5+2n)^{\prime \prime} \times (4+2m)$
    \item $(4+2n)^{\prime \prime} \times (5+2m)$
    \item $(6+2n)^{\prime \prime} \times (6+2m)$
    \item $(8+2n)^{\prime \prime} \times (4+2m)$
\end{enumerate}

We now show why the following boards are not fault-free tileable.
\begin{itemize}
    \item $\text{odd}^{\prime \prime} \times \text{odd}$
    \item $(2n)^{\prime \prime} \times 1$
    \item $(1+2n)^{\prime \prime} \times 2$
    \item $(2n)^{\prime \prime} \times 2$
    \item $1^{\prime \prime} \times (2n)$
    \item $2^{\prime \prime} \times (2+n)$
    \item $3^{\prime \prime} \times (3+2n)$
    \item $4^{\prime \prime} \times (4+2n)$
    \item $6^{\prime \prime} \times 4$
\end{itemize}

The odd by odd boards are not fault-free tileable, as before, by parity.

The $(2n)^{\prime \prime} \times 1$ boards are not fault-free tileable. Any first tile placement will result in a fault-line.

\begin{figure}[h]
     \centering
         \includegraphics[scale=.25]{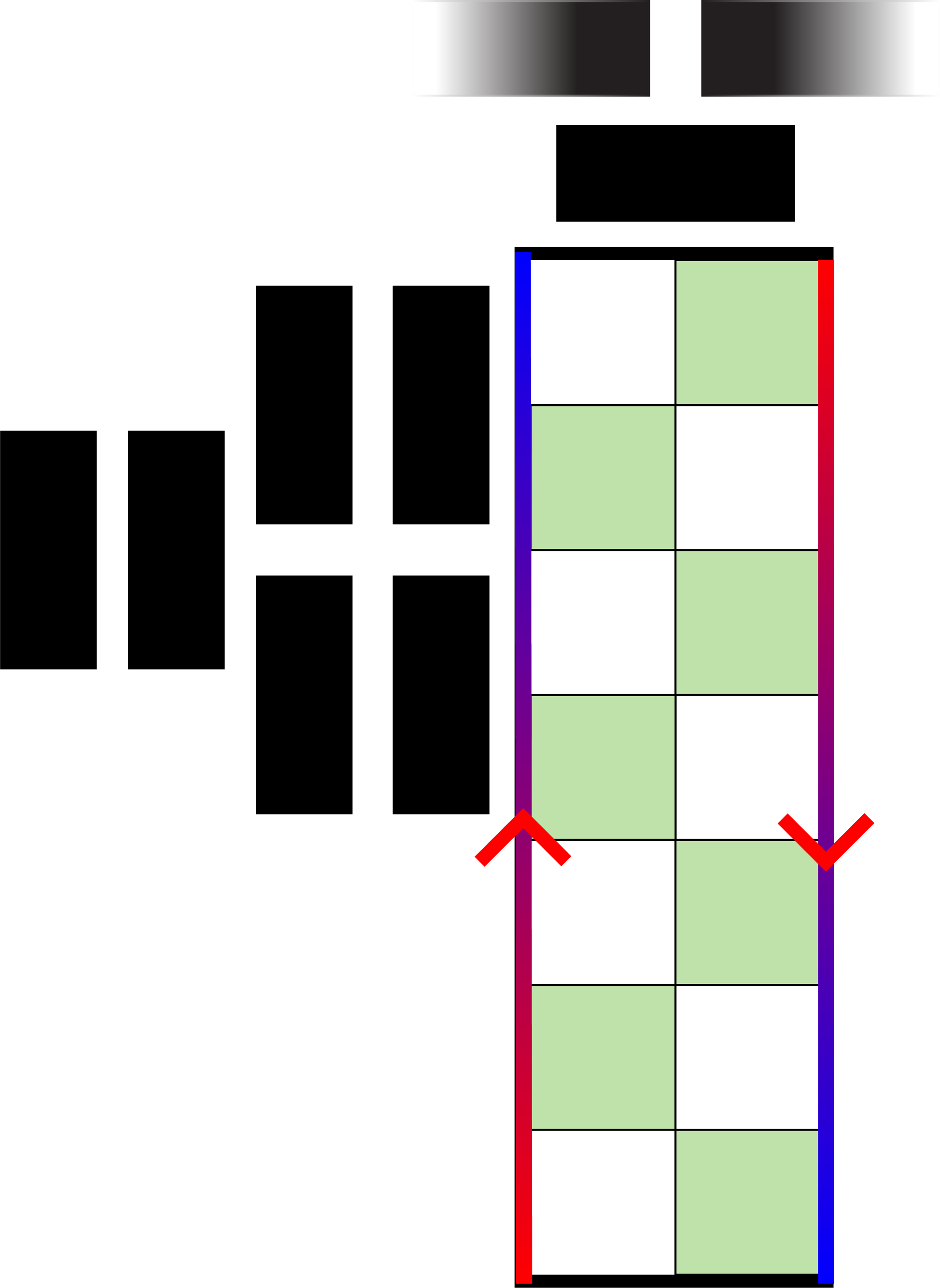}
         \caption{The $7^{\prime \prime} \times 2$ board and its forced tile sequence}
         \label{fig:7''x2}
\end{figure}

The $(1+2n)^{\prime \prime} \times 2$ boards are not fault-free tileable. An example of this is shown in Figure~\ref{fig:7''x2}. These boards require $2$ tiles to cover all vertical fault-lines and $2n$ tiles to cover all horizontal fault-lines. This is the minimum number of required tiles to maintain an even parity of remaining spaces to be covered in all rows and columns. Note the top and bottom row are the same row. This is a total of $2n+2$ tiles required. This board, however, only has room for $2n+1$ tiles. Thus, it is not fault-free tileable.

\begin{figure}[h]
     \centering
         \includegraphics[scale=.25]{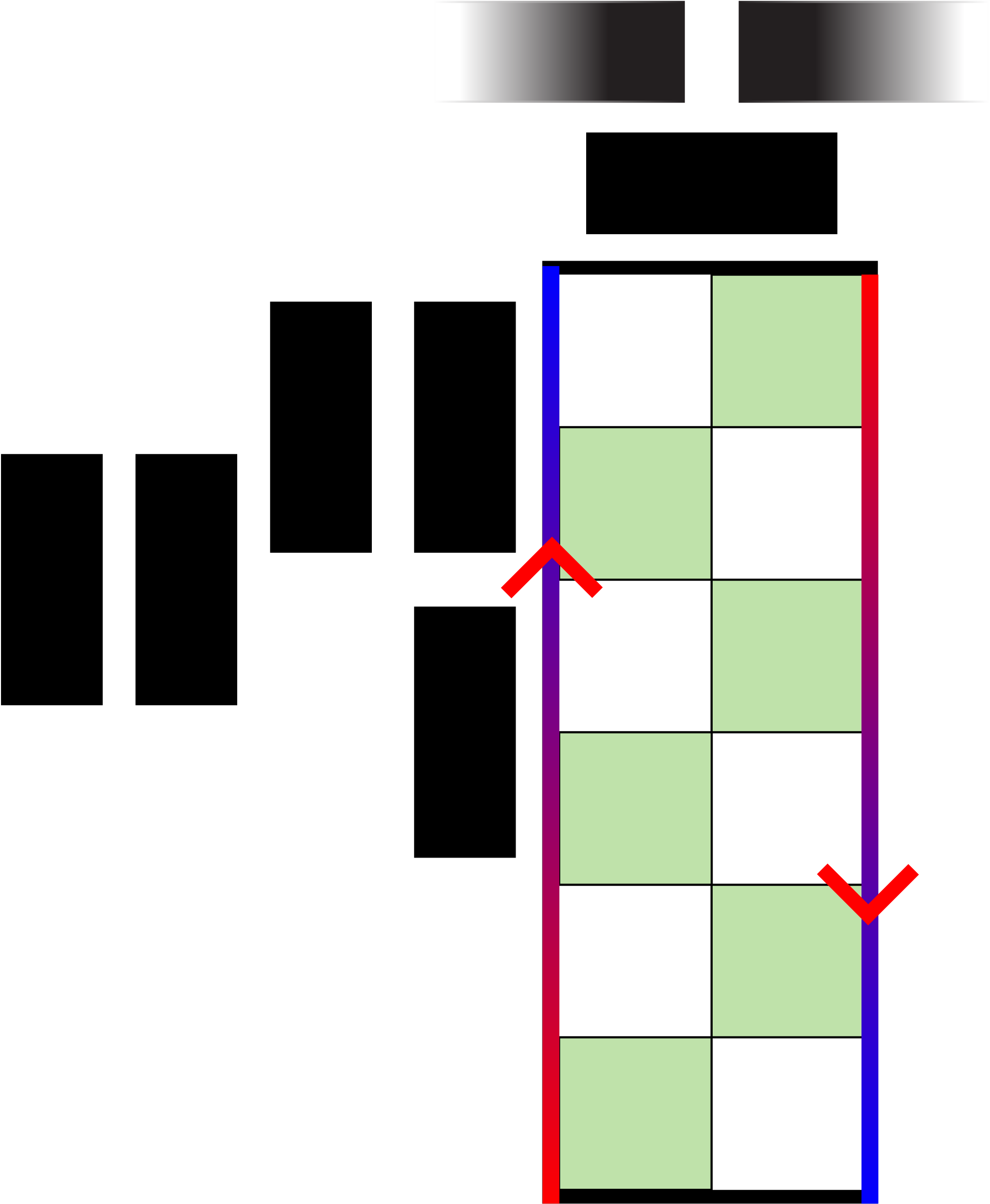}
         \caption{The $6^{\prime \prime} \times 2$ board with and its required tiles to be fault-free tileable}
         \label{fig:6''x2}
\end{figure}

Similarly, the $(2n)^{\prime \prime} \times 2$ boards are not fault-free tileable. An example on the $6^{\prime \prime} \times 2$ board is shown in Figure~\ref{fig:6''x2}. These boards require $2$ tiles to cover all vertical fault lines and $2n-1$ tiles to cover all horizontal fault-lines. This is a total of $2n+1$ tiles, while the board only has space for $2n$ tiles. Interestingly, these boards actually force the placement of more than $2n+1$ tiles. The wrapping horizontal tile will cover two spaces of the same color. This will leave an unequal number of light and dark tiles remaining. All other tiles cover one light and one dark space. Thus an additional tile must be placed to cover the two remaining spaces of the color remaining. This tile must also be a wrapping horizontal tile in order to cover $2$ tiles of the same color.

The $1^{\prime \prime} \times (2n)$ boards are not fault-free tileable. Any first tile placement will result in a fault-line.

The $2^{\prime \prime} \times (2+n)$  boards are not fault-free tileable. A vertical tile must be placed to cover the horizontal fault-line, which will result in a vertical fault-line on either side of the tile.

The $3^{\prime \prime} \times (3+2n)$ boards are not fault-free tileable. A vertical tile must be placed to cover a horizontal fault-line. This forces a horizontal tile to be placed to cover the single remaining space in the column, which will result in a fault-line.

\begin{figure}[h]
     \centering
         \includegraphics[scale=.25]{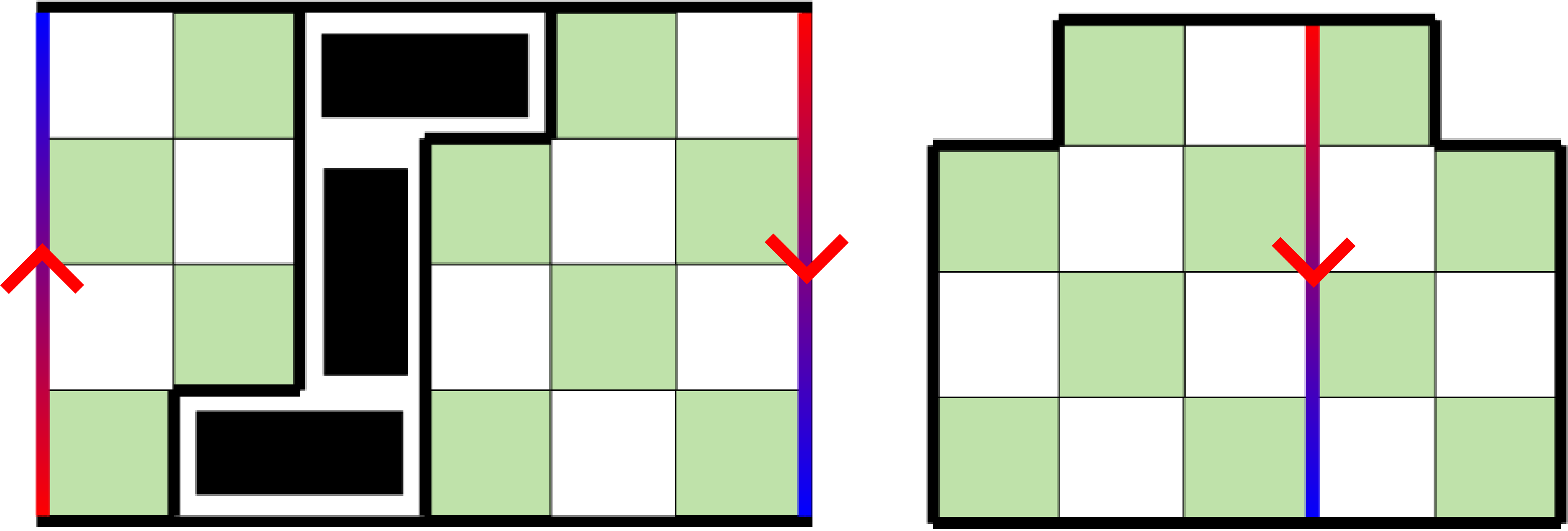}
         \caption{The $4^{\prime \prime} \times 6$ board and its forced sequence with remaining shape}
         \label{fig:6x4''ex1}
\end{figure}

The $4^{\prime \prime} \times (4+2n)$ boards are not fault-free tileable. A vertical tile must be placed to cover the middle horizontal fault-line. To avoid a vertical fault-line, we must place two tiles, one crossing each of the vertical fault-lines on opposing sides of the original tile, as shown in Figure~\ref{fig:6x4''ex1}. If we isolate the shape remaining and color it like a checker board, we see there is not an equal number of dark and light tiles. Each tile must cover 1 light tile, and 1 dark tile. Thus, this remaining shape cannot be tiled, and the original board cannot be fault-free tiled.

\begin{figure}[h]
     \centering
         \includegraphics[scale=.19]{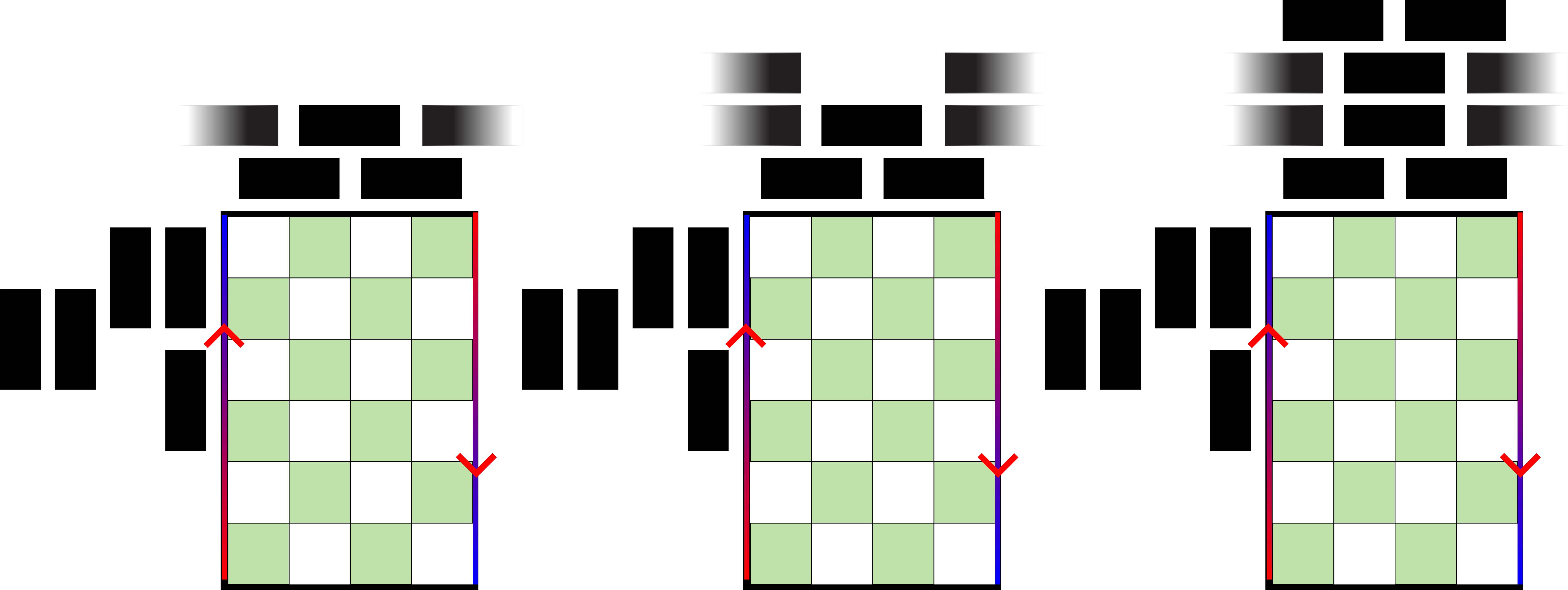}
         \caption{The $6^{\prime \prime} \times 4$ board and its required tiles to be fault-free tileable}
         \label{fig:6''x4}
\end{figure}

The $6^{\prime \prime} \times 4$ board is not fault-free tileable. This is shown in Figure~\ref{fig:6''x4}. To cover each fault line this board requires $4$ horizontal tiles and $5$ vertical tiles. The wrapping horizontal tile will cover two spaces of the same color. Every other tile will cover one light and one dark space. In order to cover the remaining two spaces of opposite color, an additional wrapping horizontal tile must be placed. To maintain parity within the columns, an additional $3$ more tiles must be added. Thus, the $6^{\prime \prime} \times 4$ board requires $13$ tiles to be fault-free tiled, while there is only space for $12$.

\begin{figure}[h]
     \centering
         \includegraphics[scale=1]{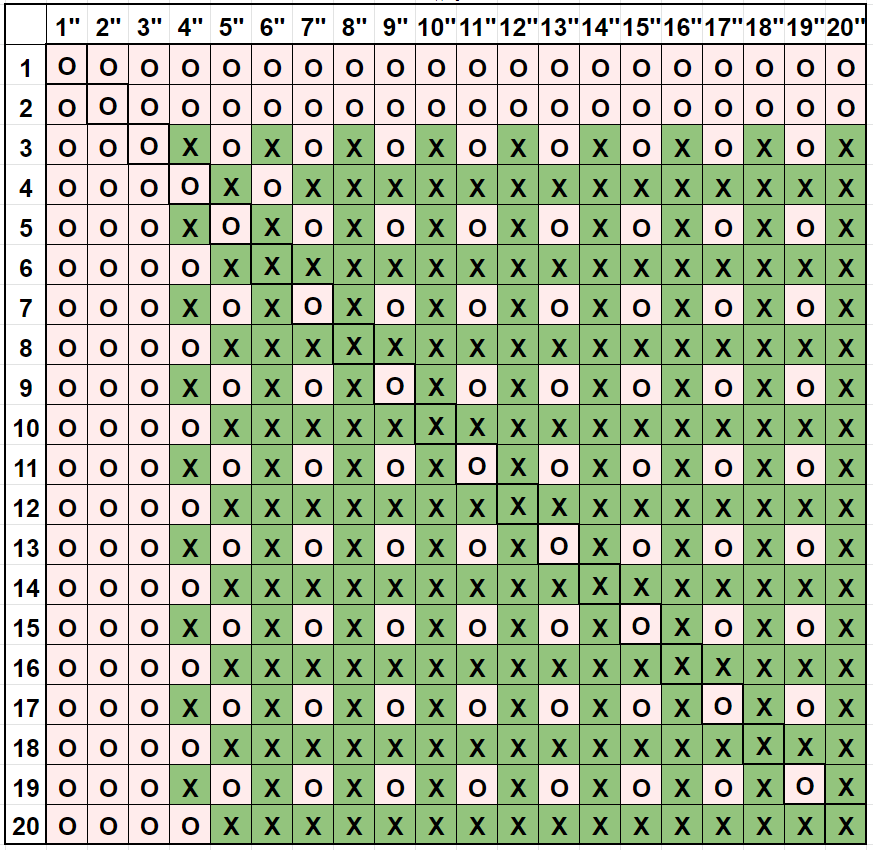}
         \caption{A chart denoting fault-free tileablity of M\"{o}bius boards where boxes marked by an X are fault-free tileable, and boxes marked by an O are not}
         \label{fig:mobiuschart}
\end{figure}

The cases given above determine the fault-free tileability of all possible M\"{o}bius strip boards. In Figure~\ref{fig:mobiuschart} the fault-free tileability of M\"{o}bius strip boards up to $20^{\prime \prime}\times 20$ are shown.

All cylinders, tori, and M\"{o}bius strips boards have been determined to be fault-free tileable or not, which was the goal of this research.

\section{Acknowledgments}
I am especially grateful to my advisor Dr. Jonathan White and to Dr. Brittney Miller for her unwavering support throughout the writing of this paper.

\newpage
\section{Figures}

\begin{figure}[h]
     \centering
         \includegraphics[width=0.9\textwidth]{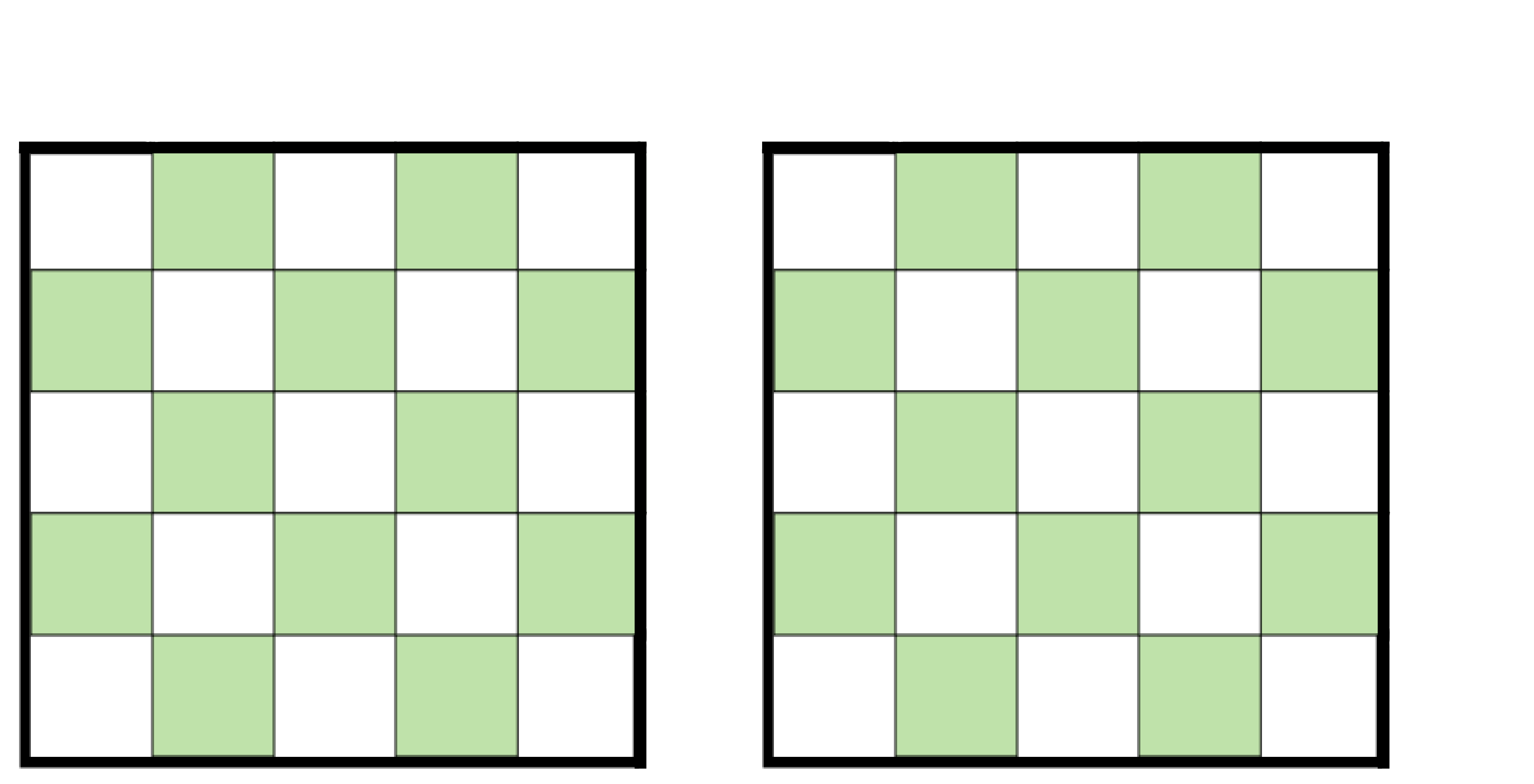}
         \caption{Blank $5\times 5$ boards for attempting tilings}
         \label{fig:5x5}
\end{figure}

\end{document}